\documentclass[journal]{IEEEtran}
\usepackage{amssymb}
\usepackage{amsfonts}
\usepackage{amsmath}
\usepackage{graphicx}
\usepackage{epsfig}
\usepackage{color}
\usepackage{subfig}
\usepackage{epstopdf}

\newtheorem{theorem}{Theorem}

\newtheorem{assumption}{Assumption}

\newtheorem{definition}[theorem]{Definition}
\newtheorem{lem}{Lemma}[section]
\newtheorem{thm}{Theorem}[section]
\newtheorem{rem}{Remark}[section]

\begin{document}
	
\title{Optimal Energy-Efficient Downlink Transmission Scheduling for
	Real-Time Wireless Networks}

\author{Lei~Miao,
	Jianfeng~Mao,
	and~Christos~G.~Cassandras,~\IEEEmembership{Fellow,~IEEE}
	\thanks{L. Miao is with the Department of Engineering Technology, Middle Tennessee State University, Murfreesboro, TN 37132 USA e-mail: lei.miao@mtsu.edu.}
	\thanks{J. Mao is with the Division of Systems and Engineering Management, Nanyang Technological University, Singapore email: jfmao@ntu.edu.sg. }
	\thanks{C. G. Cassandras is with the Division of Systems Engineering and the Deptartment of Electrical and Computer Engineering, Boston University, Brookline, MA 02446 USA email: cgc@bu.edu. }
	\thanks{The authors' work is supported in part by the National Science Foundation
		under Grant DMI-0330171, by AFOSR under grants FA9550-04-1-0133 and
		FA9550-04-1-0208, by ARO under grant DAAD19-01-0610, and by Honeywell
		Laboratories.}
	}

\maketitle
	
\begin{abstract}
		It has been shown that using appropriate channel coding schemes in wireless
		environments, transmission energy can be significantly reduced by
		controlling the packet transmission rate. This paper seeks optimal solutions
		for downlink transmission control problems, motivated by this observation
		and by the need to minimize energy consumption in real-time wireless
		networks. Our problem formulation deals with a more general setting than the
		paper authored by Gamal et. al., in which the MoveRight\ algorithm is
		proposed. The MoveRight algorithm is an iterative algorithm that converges
		to the optimal solution. We show that even under the more general setting,
		the optimal solution can be efficiently obtained through an approach
		decomposing the optimal sample path through certain \textquotedblleft
		critical tasks\textquotedblright\ which in turn can be efficiently
		identified. We include simulation results showing that our algorithm is
		significantly faster than the MoveRight algorithm. We also discuss how to
		utilize our results and receding horizon control to perform on-line
		transmission scheduling where future task information is unknown.
\end{abstract}

\begin{IEEEkeywords}
	optimization, wireless networks, energy-efficiency, real-time systems, receding horizon control.
\end{IEEEkeywords}

\IEEEpeerreviewmaketitle

\section{Introduction}

Because wireless nodes are normally powered by batteries and are expected to
remain in operation for extended periods of time, how to conserve energy in
order to extend node lifetime and network lifetime is a major research issue
in most wireless networks. One way of saving energy is to operate these
nodes at low power as long as possible. However, this will also
significantly downgrade their functionality. Therefore, there is a trade-off
between energy and the \textquotedblleft quality\textquotedblright\
delivered by wireless nodes. When \textquotedblleft
quality\textquotedblright\ is measured in terms of latency, the trade-off is
between energy and time. Examples arise in real-time computing, where a
processor trades off processing rate for energy \cite{dvs-processor}; and in
wireless transmission, where a transmitter trades off transmission speed for
energy \cite{GamNaPraUyZaInf02}.

When the energy of a wireless node is consumed mostly by communication
tasks, scheduling a RF transmission efficiently becomes extremely important
in conserving the energy of the node. It is well known that there exists an
explicit relationship between transmission power and channel capacity \cite%
{IEEEexample:Shannon}; transmission power can be adjusted by changing the
transmission rate, provided that appropriate coding schemes are used. This
provides an option to conserve the transmission energy of a wireless node by
slowing down the transmission rate. Increased latency is a direct side
effect caused by the low transmission rate and it can affect other
Quality-of-Service (QoS) metrics as well. For example, excessive delay may
cause buffer overflow, which increases the packet dropping rate. The
existence of this trade-off between energy and latency motivates \emph{%
	Dynamic Transmission Control} techniques for designing energy-efficient
wireless systems. 

To the best of our knowledge, the earliest work that captures the trade-off
between energy and latency in transmission scheduling is \cite{CollinsCruz99}%
, in which Collins and Cruz formulated a Markov decision problem for
minimizing transmission cost subject to some power constraints. By assuming
a linear dependency between transmission cost and time, their model did not
consider the potential of more energy saving by varying the transmission
rate. Berry \cite{BerryThesis} considered a Markov decision process in the
context of wireless fading channels to minimize the weighted sum of average
transmission power and a buffer cost, which corresponds to either average
delay or probability of buffer overflow. Using dynamic programming and
assuming the transmission cost to be a convex function of time, Berry
discovered some structural properties of the optimal adaptive control
policy, which relies on information on the arrival state, the queue state,
and the channel state. In \cite{AtaStaticChannel} and \cite{AtaFadingChannel}%
, Ata developed optimal dynamic power control policies subject to a QoS
constraint for Markovian queues in wireless static channels and fading
channels respectively. In his work, the optimization problem was formulated
to minimize the long-term average transmission power, given a constraint of
buffer overflow probability in equilibrium; dynamic programming and
Lagrangian relaxation approaches were used in deriving the optimal policies,
which can be expressed as functions of the packet queue length and the
channel state. Neely utilized a Lyapunov drift technique in \cite{Neely06TIT}
to develop a dynamic power allocation and routing algorithm that minimizes
the average power of a cell-partitioned wireless network. It was shown that
the on-line algorithm operates without knowledge of traffic rates or channel
statistics, and yields average power that is arbitrarily close to the
off-line optimal solution. A related problem of maximizing throughput
subject to peak and average power constraints was also discussed in \cite%
{Neely06TIT}.

Today's real-time data communications require Quality-of-Service (QoS)
guarantee for each individual packet. Another line of research aims at
minimizing the transmission energy over a single wireless link while
providing QoS guarantee. In particular, it is assumed that each packet is
associated with an arrival time (generally random), a number of bits, a hard
deadline that must be met, and an energy function. This line of work was
initially studied in \cite{TON-Energy-Efficient_Wireless} with follow-up
work in \cite{GamNaPraUyZaInf02} where a "homogeneous" case is considered
assuming all packets have the same deadline and number of bits. By
identifying some properties of this convex optimization problem, Gamal et
al. proposed the "MoveRight" algorithm in \cite{GamNaPraUyZaInf02} to solve
it iteratively. However, the rate of convergence of the MoveRight algorithm
is only obtainable for a special case of the problem when all packets have
identical energy functions; in general the MoveRight algorithm may converge
slowly. Zafer et al. \cite{Zafer09TAC} studied an optimal rate control
problem over a time-varying wireless channel, in which the channel state was
considered to be a Markov process. In particular, they considered the
scenario that $B$ units of data must be transmitted by a common deadline $T,$
and they obtained an optimal rate-control policy that minimizes the total
energy expenditure subject to short-term average power constraints. In \cite%
{Zafer07ITA} and \cite{Zafer08TIT}, the case of identical arrival time and
individual deadline is studied by Zafer et. al. In \cite{NeedlyInfocom07},
the case of identical packet size and identical delay constraint is studied
by Neely et. al. They extended the result for the case of individual packet
size and identical delay constraint in \cite{NeedlyWirelessNetworks09}. In 
\cite{ZaferTON09}, Zafer et. al. used a graphical approach to analyze the
case that each packet has its own arrival time and deadline. However, there
were certain restrictions in their setting, for example, the packet that
arrives later must have later deadlines. 
\color{black}
Wang and Li \cite{WangToWC2013} analyzed scheduling problems for bursty
packets with strict deadlines over a single time-varying wireless channel.
Assuming slotted transmission and changeable packet transmission order, they
are able to exploit structural properties of the problem to come up with an
algorithm that solves the off-line problem. In \cite{PoulakisTOVT2013},
Poulakis et. al. also studied energy efficient scheduling problems for a
single time-varying wireless channel. They considered a finite-horizon
problem where each packet must be transmitted before $D_{\max }.$ Optimal
stopping theory was used to find the optimal start transmission time between 
$[0,$ $D_{\max }]$ so as to minimize the expected energy consumption and the
average energy consumption per unit of time. In \cite{ZouToWC2013}, an
energy-efficient and deadline-constrained problem was formulated in lossy
networks to maximize the probability that a packet is delivered within the
deadline minus a transmission energy cost. Dynamic programming based solutions were developed under a
finite-state Markov channel model. Shan et. al. \cite{ShanJSAC2015} studied
discrete rate scheduling problems for packets with individual deadlines in
energy harvesting systems. Under the assumption that later packet arrivals
have later deadlines, they established connections between continuous rate
and discrete rate algorithms. A truncation algorithm was also developed to
handle the case that harvested energy is insufficient to guarantee all
packets' deadlines are met. Tomasi et. al. \cite{TomasiToWC2015} developed
transmission strategies to deliver a prescribed number of packets by a
common deadline $T$ while minimizing transmission attempts. Modeling the
time-varying correlated wireless channel as a Markov chain, they used
dynamic programming and a heuristic strategy to address three systems, in
which the receiver provides the channel state information to the transmitter
differently. Zhong and Xu \cite{ZhongXuInfocom2008} formulated optimization problems that minimize the energy consumption of a set of tasks with task-dependent energy functions and packet lengths. In their problem formulation, the energy functions include both transmission energy and circuit power consumption. To obtain the optimal solution for the off-line case with backlogged tasks only, they developed an iterative algorithm RADB whose complexity is $O(n^{2})$ ($n$ is the number of tasks). The authors show via simulation that the RADB algorithm achieves good performance when used in on-line scheduling. In \cite{VazeInfocom2013}, Vaze derived the competitive ratios of on-line transmission scheduling algorithms for single-source and two-source Gaussian channels in energy harvesting systems. In Vaze's problem formulation, the goal is to minimize the transmission time of fixed $B$ bits using harvested energy, which arrive in chunks randomly.

In the above papers, the closest ones to this paper are \cite%
{GamNaPraUyZaInf02}, \cite{NeedlyWirelessNetworks09}, and \cite{ZaferTON09}.
In this paper, we consider the transmission control problem in the scenario
that each task has arbitrary arrival time, deadline, and number of bits.
Therefore, the problem we study in this paper is more generic and
challenging. 

Our model also allows each packet to have its own energy function. This
makes our results especially applicable to Download Transmission Scheduling
(DTS) scenarios, where a transmitter transmits to multiple receivers over
slow-fading channels. Our contributions are the following: by analyzing the
structure of the optimal sample path, we solve the DTS problem efficiently
using a two-fold decomposition approach. First, we establish that the
problem can be reduced to a set of subproblems over segments of the optimal
sample path defined by \textquotedblleft critical tasks\textquotedblright .
Secondly, we establish that solving each subproblem boils down to solving
nonlinear algebraic equations for the corresponding segments. Based on the
above decomposition approach, an efficient algorithm that solves the DTS
problem is proposed and compared to the MoveRight algorithm. Simulation results show that our algorithm is typically an order of
magnitude faster than the MoveRight algorithm.

The main results of the paper were previously published in \cite%
{MiaoCassInfocom05}. In this journal version, we have improved most proofs
and moved them to an appendix in order to enhance the continuity of the
analysis in the paper. We have added summaries and explanations between the
technical results to enhance its readability. 

We have added Section III.C, in which the maximum power constraints are
added and discussed. 

In addition, we have added Section IV. In this section, we discuss how to
use our algorithm and Receding Horizon Control to perform on-line
transmission scheduling where the task information is unknown. New
simulation results are also provided in this section.

The structure of the paper is the following: in Section II, we formulate our
DTS\ problem and discuss some related work; the main results of DTS are
presented in Section III, where an efficient algorithm is proposed and shown
to be optimal; in Section IV, we discuss how our main results can be used to
perform on-line transmission control; finally, we conclude in Section V.

\section{The Downlink Transmission Scheduling Problem and Related Work}

We assume the channel between the transmitter and the receiver is an
Additive White Gaussian Noise (AWGN) channel and the interference to the
receiver is negligible. The received signal at time $t$ can be written as:%
\begin{equation}
	Y(t)=\sqrt{g(t)}X(t)+n(t),  \label{received_signal}
\end{equation}%
where $g(t)$ is the channel gain, $X(t)$ is the transmitted signal, and $%
n(t) $ is additive white Gaussian noise \cite{Goldsmith97}. Note that in
this section we consider the case when the transmitter is in isolation from
other transmitters so that the interference is negligible. Due to channel
fading, $g(t)$ is time-varying in general. We will consider $g(t)$ to be
time-invariant during the transmission of a single packet. Although in
practice the channel state may change during the transmission of a packet,
our results are still helpful, since, it is valid to estimate unknown future
channel state to be static for each packet in an on-line setting. Note that
our results can be possibly extended to fast fading channels as well.

The DTS problem arises when a wireless node has a set of \emph{N} packets
that need to be sent to different neighboring nodes. The goal is to minimize
the total transmission energy consumption while guaranteeing hard deadline
satisfaction for each individual packet. Since each packet can be considered
as a communication task, we the terms \textquotedblleft
task\textquotedblright\ and \textquotedblleft packet\textquotedblright\
interchangeably in what follows. We model the transmitter as a single-server
queueing system operating on a nonpreemptive and First-Come-First-Ferved
(FCFS) basis, whose dynamics are given by the well-known max-plus equation%
\begin{equation}
	x_{i}=\max (x_{i-1},a_{i})+s_{i}  \label{max_plus_equation}
\end{equation}%
where $a_{i}$ is the arrival time of task $i=1,2,\ldots ,$ $x_{i}$ is the
time when task $i$ completes service, and $s_{i}$ is its (generally random)
service time.

Note that although preemption is often easy and straightforward in computing
systems, it is very costly and also technically hard in wireless
transmissions. Therefore, we assume a nonpreemptive model in this paper.
Transmission rate control typically occurs in the physical layer, and
changing packet order may cause problems in the upper layers of the network
stack. Thus, we use a simple FCFS model to avoid packet out-of-sequence
problems. It is also worth noting that even if the packet order is
changeable, determining the optimal packet order is a separate problem. Once
the order of transmission is decided by a specific scheduling policy, our
work can be used to minimize the energy expenditure for that specific order.

The service time $s_{i}$ is controlled by the transmission rate, which is
determined by transmission power and coding scheme. However, it turns out
that it is more convenient to use the reciprocal of the transmission rate as
our control variable in the DTS\ problem. Thus, we define $\tau $ to be the
transmission time per bit and $\omega _{i}(\tau )$ to be the energy cost per
bit for task $i$. Clearly, $\omega _{i}(\tau )$ is a function of $\tau $.
Since the channel gain $g(t)$ in (\ref{received_signal}) is constant, $%
\omega _{i}(\tau )$ is kept fixed during the transmission of task $i$.

We formulate the off-line DTS problem as follows:%
\begin{equation*}
	\begin{tabular}{ll}
		\textbf{P1:} & $\underset{\tau _{1},\ldots ,\tau _{N}}{\min }%
		\sum_{i=1}^{N}v_{i}\omega _{i}(\tau _{i})$ \\ 
		$s.t.$ & $x_{i}=\max (x_{i-1},a_{i})+v_{i}\tau _{i}\leq d_{i},$ $i=1,\ldots
		N $ \\ 
		& $\tau _{i}>0,\text{ }x_{0}=0.$%
	\end{tabular}%
\end{equation*}%
where $d_{i}$ and $v_{i}$ are the deadline and the number of bits of task $i$
respectively.

In realistic scenarios, the maximum transmission power of a wireless system
puts a constraint on each $\tau _{i},$ i.e., $\tau _{i}\geq \tau _{i\_\min
}, $ where $\tau _{i\_\min }$ is the minimum amount of time used for
transmitting one bit in task $i$. For ease of analysis, we omitted this
constraint in \textbf{P1. }However, it is important to note that special
handling is needed in real-world systems for the case that the optimal
solution $\tau _{i}^{\ast }$ is below the minimum value $\tau _{i\_\min }$.
For example, the system may simply choose to drop the packet or transmit the
packet using control $\tau _{i\_\min }.$ We will discuss the problem that
includes this constraint in Section III.C and Section IV. 

Note that in the off-line setting, we consider $a_{i},$ $d_{i}$ and $v_{i}$
are known. The downlink scheduling problem formulated in \cite%
{GamNaPraUyZaInf02} is a special case of \textbf{P1} above: in \cite%
{GamNaPraUyZaInf02} each task has the same deadline and number of bits,
i.e., $d_{i}=T,v_{i}=v,$ for all $i$. Note that transmission rate
constraints are omitted in \textbf{P1} and we assume the transmission rate
can vary continuously. In practical systems, the control can always be
rounded to the nearest achievable value \cite{KomSno03}.

Problem \textbf{P1} above is similar to the general class of problems
studied in \cite{cgc-forward} and \cite{cgc-forward-improved} without the
constraints $x_{i}\leq d_{i}$, where a decomposition algorithm termed the
Forward Algorithm (FA) was derived. As shown in \cite{cgc-forward} and \cite%
{cgc-forward-improved}, instead of solving this complex nonlinear
optimization problem, we can decompose the optimal sample path into a number
of busy periods. A \emph{busy period} (BP) is a contiguous set of tasks $%
\{k,...,n\}$ such that the following three conditions are satisfied: $%
x_{k-1}<a_{k}$, $x_{n}<a_{n+1}$, and $x_{i}\geq a_{i+1}$, for every $%
i=k,\ldots ,n-1$. Notice that \textbf{P1} above
	exploits static control ($\tau _{i}$\ kept fixed
	during the service time of task $i$). This is
	straightforward in wireless transmission control since the transmission rate
	of a single packet/task is often fixed. In addition, it has been shown in 
	\cite{MiaoCas2004TN} that when the energy functions $\omega _{i}(\tau ),$%
\ $i=1,\ldots N,$\ are strictly	convex and monotonically decreasing in $\tau ,$\
	there is no benefit in applying dynamic control ($\tau _{i}$%
\ varies over time during the service time of task $i$%
). It has also been shown in \cite{MaoCTDA} that when
the energy functions are identical in \textbf{P1}, its solution is obtained
by an efficient algorithm (Critical Task Decomposition Algorithm) that
decomposes the optimal sample path even further and does not require solving
any convex optimization problem at all. In this paper, we will consider the
much harder case that the energy functions are \emph{task-dependent}. When
the energy functions are homogeneous, it is shown in \cite{MaoCTDA} that the
exact form of the energy function does not matter in finding the optimal
solutions. The main challenge of having heterogeneous energy functions is
that these energy functions will be used to identify the optimal solutions,
and this adds an extra layer of complexity. We shall still use the
decomposition idea in \cite{MaoCTDA}, and we will use $\{\tau _{i}^{\ast }\}$
and $\{x_{i}^{\ast }\}$, $i=1,\ldots ,N$, to denote the optimal solution of 
\textbf{P1 }and the corresponding task departure times respectively.

Typically, $\omega _{i}(\tau )$ is determined by factors including the
channel gain $g(t)$, transmission distance, signal to noise ratio, and so
on. Therefore, when a wireless node transmits to different neighbors at
different time, different $\omega _{i}(\tau )$ are involved. We begin with
an assumption that will be made throughout our analysis.

\begin{assumption}
	\label{A2}In AWGN channels, $\omega _{i}(\tau )$ is nonnegative, strictly
	convex, monotonically decreasing, differentiable, and $\lim_{\tau
		\rightarrow 0}\dot{\omega}_{i}(\tau )=$ $-\infty $ .
\end{assumption}

Assumption \ref{A2} is justified in \cite{GamNaPraUyZaInf02} and channel
coding schemes supporting this assumption can be found in \cite%
{TON-Energy-Efficient_Wireless}. Note that the result obtained in \cite%
{MiaoCas2004TN} can be readily applied here: the unique optimal control to 
\textbf{P1} is static. This means that we do not need to vary the
transmission rate of task $i$ during its transmission time.

\section{Main Results of DTS\label{main_results}}

\subsection{Optimal Sample Path Decomposition}

The following two lemmas help us to decompose the optimal sample path of 
\textbf{P1}. Their proofs are very similar to the proofs for Lemmas 1 in 
\cite{MaoCTDA}, and only monotonicity of $\omega _{i}(\tau )$ is required.
We omit the proofs here.

\begin{lem}
	\label{DepartureIsDeadline}If $d_{i}<a_{i+1},$ then $x_{i}^{\ast}=d_{i}.$
\end{lem}

\begin{lem}
	\label{DepartureLargerThanNextArrival}If $d_{i}\geq a_{i+1},$ then $%
	a_{i+1}\leq x_{i}^{\ast}.$
\end{lem}

Recalling the definition of a BP, Lemmas \ref{DepartureIsDeadline}, \ref%
{DepartureLargerThanNextArrival} show that the BP optimal structure can be
explicitly determined by the deadline-arrival relationship, i.e., a sequence
of contiguous packets $\{k,\ldots ,n\}$ is a BP if and only if the following
is satisfied: $d_{k-1}<a_{k},$ $d_{n}<a_{n+1},$ $d_{i}\geq a_{i+1},$ for all 
$i\in \{k,\ldots ,n-1\}.$ After identifying each BP on the optimal sample
path, problem \textbf{P1} is reduced to solving a separate problem over each
BP. We formulate the following optimization problem for BP $\{k,\ldots ,n\}.$%
\begin{equation*}
	\begin{tabular}{ll}
		$Q(k,n):$ & $\underset{\tau _{k},\ldots ,\tau _{n}}{\min }%
		\sum_{i=k}^{n}v_{i}\omega _{i}(\tau _{i})$ \\ 
		$\ \ \ s.t.$ & $x_{i}=a_{k}+\sum_{j=k}^{i}v_{i}\tau _{i}\leq d_{i},$ $%
		i=k,\ldots ,n$ \\ 
		& 
		$\tau _{i}>0,$ $i=k,\ldots ,n,$%
		\\ 
		& $x_{i}\geq a_{i+1},\text{ }i=k,\ldots ,n-1.$%
	\end{tabular}%
\end{equation*}

Although $Q(k,n)$ is easier than \textbf{P1} (since it does not contain
max-plus equations, which are nondifferentiable), it is still a hard convex
optimization problem. Naturally, we would like to solve $Q(k,n)$
efficiently. As we will show, this is indeed possible by further decomposing
a BP $\{k,\ldots,n\}$ through special tasks called ``critical tasks'', which
are defined as follows:

\begin{definition}
	\label{Definition1}Suppose both task $i$ and $i+1$ are within a BP $%
	\{k,\ldots,n\}$ on the optimal sample path of \textbf{P1}. If $\dot{\omega }%
	_{i}(\tau_{i}^{\ast})\neq\dot{\omega}_{i+1}(\tau_{i+1}^{\ast}),$ task $i $
	is \textbf{critical}. If $\dot{\omega}_{i}(\tau_{i}^{\ast})>\dot{\omega}%
	_{i+1}(\tau_{i+1}^{\ast}),$ then task $i$ is \textbf{left-critical}. If $%
	\dot{\omega}_{i}(\tau_{i}^{\ast})<\dot{\omega}_{i+1}(\tau_{i+1}^{\ast}),$
	then task $i$ is \textbf{right-critical}.
\end{definition}

These critical tasks are special because the derivatives of the energy
function change after these tasks are transmitted on the optimal sample
path. Therefore, identifying critical tasks is crucial in solving $Q(k,n).$
In fact, Gamal et al. \cite{GamNaPraUyZaInf02} observed the existence of
left-critical tasks. However, they did not make use of them in
characterizing the optimal sample path. In order to accomplish this, we need
to study the relationship between critical tasks and the structure of the
optimal sample path. An auxiliary lemma will be introduced first.

\begin{lem}
	\label{Convex_Transmission}If $v_{1}\tau _{1}+v_{2}\tau _{2}=v_{1}\tau
	_{1}^{^{\prime }}+v_{2}\tau _{2}^{^{\prime }},$ $\tau _{1}^{^{\prime }}<\tau
	_{1},\tau _{2}^{^{\prime }}>\tau _{2},$ and $\dot{\omega}_{1}(\tau
	_{1}^{^{\prime }})>\dot{\omega}_{2}(\tau _{2}^{^{\prime }}),$ then, $%
	v_{1}\omega _{1}(\tau _{1})+v_{2}\omega _{2}(\tau _{2})>v_{1}\omega
	_{1}(\tau _{1}^{^{\prime }})+v_{2}\omega _{2}(\tau _{2}^{^{\prime }}).$
\end{lem}

Lemma \ref{Convex_Transmission} implies that under Assumption \ref{A2}
(especially, the convexity assumption), it takes the least amount of energy
to transmit two tasks in a given amount of time when the derivatives of the
two energy functions have the least amount of difference. As we will see
later, this auxiliary lemma will be used to establish other important
results. Next, we will discuss what exactly makes the critical tasks
(defined in Definition \ref{Definition1}) special.

\begin{lem}
	\label{Critical_Condition}Suppose both task $i$ and $i+1$ are within a BP $%
	\{k,\ldots ,n\}$ on the optimal sample path of \textbf{P1}$.$ \emph{(i)} If
	task $i$ is left-critical, then $x_{i}^{\ast }=a_{i+1}.$ \emph{(ii)} If task 
	$i$ is right-critical, then $x_{i}^{\ast }=d_{i}.$
\end{lem}

This result shows that if a task is \emph{left-critical} or \emph{%
	right-critical} on the optimal sample path, its optimal departure time is
given by the next arrival time or its deadline respectively. The lemma
implies that when $a_{i+1}<x_{i}^{\ast }<d_{i},$ task $i$ is neither \emph{%
	left-critical} nor \emph{right-critical}. In our next step, we will study
the commonality among a block of consecutive non-critical tasks, which are
in the middle of two adjacent critical tasks. By hoping so, we will have a
better understanding of the structure of the optimal sample path, using
which we will develop an efficient algorithm to solve $Q(k,n).$

\begin{rem}
	\label{corollary_rate_equal}For any two neighboring tasks $i$ and $i+1$ in a
	BP $\{k,\ldots ,$ $n\}$ on the optimal sample path of \textbf{P1}, if task $%
	i $ is not a critical task, then $\dot{\omega}_{i}(\tau _{i}^{\ast })=\dot{%
		\omega}_{i+1}(\tau _{i+1}^{\ast }).$
\end{rem}

\begin{figure}[tbp]
\begin{center}  
\includegraphics[height=2.0951in,width=3in,angle=0]{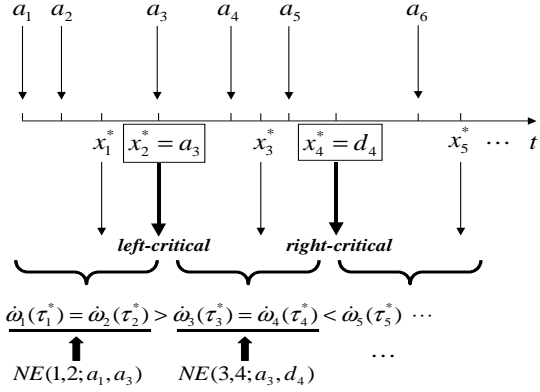}  
\caption{\small \sl Illustration of the optimal
	structure of BP \{k,\ldots ,n\}. \label{Fig_Optimal_Structure}}  
\end{center}  
\end{figure}

This remark\textbf{\ }is the direct result of Definition \ref{Definition1}.
Using this remark and Lemma \ref{Critical_Condition}, \ we can obtain the
structure of BP $\{k,\ldots ,n\}$ on the optimal sample path of \textbf{P1 }%
as follows: $\{k,\ldots ,n\}$ is characterized by a sequence of tasks $%
S=\{c_{0},\ldots ,c_{m+1}\},$ in which $c_{0}=k$, $\{c_{1},\ldots ,c_{m}\}$
contains all critical tasks in $\{k,\ldots ,n\}$ (the optimal departure
times of these critical tasks are given by Lemma \ref{Critical_Condition}),
and task $c_{m+1}=n$. Moreover, let $c_{i}$, $c_{i+1}$ be adjacent tasks in $%
S$. Then, the segment of tasks 
\begin{equation}
	\left\{ 
	\begin{tabular}{ll}
		$\{c_{i},\ldots ,c_{i+1}\},$ & if $i=0$ \\ 
		$\{c_{i}+1,\ldots ,c_{i+1}\},$ & if $0<i\leq m$%
	\end{tabular}%
	\ \right.  \label{segment_definition}
\end{equation}%
is operated at some $\tau $ such that the derivatives of their energy
functions are all the same. To have a better understanding of this optimal
structure, see Fig. \ref{Fig_Optimal_Structure}. In this example, task $2$
is left-critical and task $4$ is right-critical. Their optimal departure
times are $a_{3}$ and $d_{4}$ respectively. \ In the set $S=\{1,2,4,...\}$,
tasks $\{1,2\}$ and $\{3,4\}$ are examples of the segments defined above.
Invoking Remark \ref{corollary_rate_equal}, $\tau _{1}^{\ast },\ldots ,\tau
_{4}^{\ast }$ are characterized by $\dot{\omega}_{1}(\tau _{1}^{\ast })=\dot{%
	\omega}_{2}(\tau _{2}^{\ast }),$ $\dot{\omega}_{3}(\tau _{3}^{\ast })=\dot{%
	\omega}_{4}(\tau _{4}^{\ast })$ and $\dot{\omega}_{2}(\tau _{2}^{\ast })>%
\dot{\omega}_{3}(\tau _{3}^{\ast }),$ $\dot{\omega}_{4}(\tau _{4}^{\ast })<%
\dot{\omega}_{5}(\tau _{5}^{\ast }).$

In order to obtain our main result of this section and the explicit algorithm that solves $Q(k,n),$ we
define next a system of nonlinear algebraic equations as follows with $i<j,$ 
$0\leq t_{1}\leq t_{2},$ and unknown variables $\tau_{i},\ldots,\tau_{j}:$%
\begin{equation*}
\begin{tabular}{ll}
$NE(i,j;t_{1},t_{2}):$ & $\sum_{m=i}^{j}\tau_{m}v_{m}=t_{2}-t_{1},$ \\ 
& $\dot{\omega}_{m}(\tau_{m})=\dot{\omega}_{m+1}(\tau_{m+1}),$ \\ 
& $m=i,\ldots,j-1.$%
\end{tabular}%
\end{equation*}

Its solution minimizes the total energy of transmitting tasks $\{i,\ldots,$ $%
j\}$ that do not contain critical tasks within time interval $t_{2}-t_{1}$.
Note that when $i=j$, the above nonlinear algebraic equations reduce to a\
single linear equation $\tau_{i}v_{i}=t_{2}-t_{1}.$

In Fig. \ref{Fig_Optimal_Structure}, we illustrated the structure of a BP on
the optimal sample path of \textbf{P1}. In fact, given all critical tasks in
the BP, the optimal solution can be obtained by solving a set of $NE$
systems, one for each segment defined in (\ref{segment_definition}). For
example, in Fig. \ref{Fig_Optimal_Structure}, the optimal controls of tasks $%
\{1,2\}$ and $\{3,4\}$ can be obtained by solving $NE(1,2;a_{1},a_{3})$ and $%
NE(3,4;a_{3},d_{4})$ respectively.

At this point, we have established that solving problem $Q(k,n)$ boils down
to identifying critical tasks on its optimal sample path. This relies on
some additional properties of the optimal sample path. To obtain them, we
need to first study the properties of $NE(i,j;t_{1},t_{2}).$

We denote the solution to $NE(i,j;t_{1},t_{2})$ by $\tau
_{i}(t_{1},t_{2}),\ldots ,\tau _{j}(t_{1},t_{2})$. We define the common
derivative in $NE(i,j;t_{1},t_{2})$: 
\begin{equation*}
\sigma _{i,j}(t_{1},t_{2})=\dot{\omega}_{m}(\tau _{m}(t_{1},t_{2})),\text{
	for any }m\text{, }i\leq m\leq j.
\end{equation*}%
and note that $\sigma _{i,j}(t_{1},t_{2})$ is the derivative of the energy
function of any task in $\{i,\ldots ,j\}$. When $t_{1}=t_{2},$ we set $%
\sigma _{i,j}(t_{1},t_{2})$ to $-\infty $. Later, when invoking the
definition of critical tasks, we will use $\sigma _{i,j}(t_{1},t_{2})$
instead of the derivative of the energy function of a single task.

Now, we are ready to introduce the properties of $NE(i,j;t_{1},t_{2})$ in
the next lemma.

\begin{lem}
\label{NE_Properties}When $t_{1}<t_{2},$ $NE(i,j;t_{1},t_{2})$ has the
following properties:

\emph{(i)} It has a unique solution.

\emph{(ii)} The common derivative $\sigma_{i,j}(t_{1},t_{2})$ is a
monotonically increasing function of $\Delta=t_{2}-t_{1},$ $i.e.,$ 
\begin{equation*}
\sigma_{i,j}(t_{1},t_{2})<\sigma_{i,j}(t_{3},t_{4}),\text{ if }%
t_{4}-t_{3}>t_{2}-t_{1}.
\end{equation*}

\emph{(iii)} For any $p$, $i\leq p<j,$ define the partial sum $%
S_{ip}\equiv\sum_{m=i}^{p}\tau_{m}(t_{1},$ $t_{2})v_{m}.$ Then, 
\begin{equation*}
\sigma_{i,p}(t_{1},t_{1}+S_{ip})=\sigma_{p+1,j}(t_{1}+S_{ip},t_{2})=%
\sigma_{i,j}(t_{1},t_{2})
\end{equation*}

\emph{(iv)} For any $p$, $i\leq p<j,$ $t_{1}<t_{3}<t_{2},$ let $c_{1}=\sigma
_{i,p}(t_{1},t_{3}),$ $c_{2}=\sigma _{p+1,j}(t_{3},t_{2}),$ $c_{3}=\sigma
_{i,j}(t_{1},t_{2}).$ If $c_{q}\neq c_{r}$ $\forall q,r\in \{1,2,3\},q\neq
r, $ then $\min (c_{1},c_{2})<c_{3}<\max (c_{1},c_{2}).$
\end{lem}

\subsection{Left and Right-critical Task Identification}

Based on the above results, we have characterized the special structure of
the optimal sample path of \textbf{P1}. To summarize, Lemmas \ref%
{DepartureIsDeadline} and \ref{DepartureLargerThanNextArrival} show that the
BP structure of the optimal sample path can be explicitly determined by the
deadline-arrival relationship. This transforms \textbf{P1} into a set of
simpler convex optimization problems with linear constraints. Although the
problem becomes easier to solve, it is still computationally hard for
wireless devices without powerful processors and sufficient energy. Note
that in the homogeneous case, when all tasks have the same arrival time and
deadline, they should be transmitted with the same derivatives of their cost
functions. In this case, the optimal solution can be obtained by solving the
nonlinear system $NE(i,j;t_{1},t_{2}).$ With the presence of inhomogeneous
real-time constraints, we showed in Lemma \ref{Critical_Condition} and
Remark \ref{corollary_rate_equal} that a set of ``critical tasks'' play a
key role to determine the optimal sample path, i.e., the derivatives of the
cost functions only change at these critical tasks. Once they are
determined, the original problem $Q(k,n)$ boils down to set of nonlinear
algebraic equations.

Having obtained the properties of $NE(i,j;t_{1},t_{2}),$ we will next
develop an efficient algorithm to identify critical tasks. Without loss of
generality, we only prove the correctness of identifying the first critical
task. Other critical tasks can be identified iteratively. In addition, our
proof will focus on right-critical tasks only, and we omit the proof for
left-critical tasks, which is very similar.

We will first give some definitions. For tasks $(p,i)$ within a BP $%
\{k,\ldots ,n\}$, i.e., $k\leq p<i\leq n,$ define:%
\begin{equation*}
T_{1}(k,p)=\left\{ 
\begin{array}{c}
a_{k},\text{ }p=k\text{ } \\ 
x_{p-1}^{\ast },\text{ }p>k%
\end{array}%
\right.
\end{equation*}%
\begin{equation*}
T_{2}(n,i)=\left\{ 
\begin{array}{c}
a_{i+1},\text{ }i<n \\ 
d_{n},\text{ }i=n%
\end{array}%
\right.
\end{equation*}

\ Recalling the definition of a BP, $T_{1}(k,p)$ is defined as \emph{the
	optimal starting transmission time for task }$p$, which is within a BP
starting with task $k$. Recalling Lemmas \ref{DepartureIsDeadline} and \ref%
{DepartureLargerThanNextArrival}, $T_{2}(n,i)$ is defined as \emph{the
	earliest possible transmission ending time for task }$i$, which is within a
BP ending with task $n$. Note that in order to guarantee the real-time
constraints, task $i$ must be done by its deadline $d_{i}.$ We will use $%
T_{1}(k,p),$ $T_{2}(n,i)$, and $d_{i}$ later to identify critical tasks.

We further define: 
\begin{gather}
R_{i}=\underset{s\in\{p,\ldots,i-1\}}{\arg\max}\{%
\begin{array}{c}
\sigma_{p,s}(T_{1}(k,p),d_{s}) \\ 
\leq\sigma_{p,j}(T_{1}(k,p),d_{j}),%
\end{array}
\label{Ri} \\
\text{for }i\text{, }p<i\leq n,\text{ and all }j\in\{p,\ldots,i-1\}\}  \notag
\end{gather}

\begin{gather}
L_{i}=\underset{s\in\{p,\ldots,i-1\}}{\arg\max}\{%
\begin{array}{c}
\text{ }\sigma_{p,s}(T_{1}(k,p),T_{2}(n,s)) \\ 
\geq\sigma_{p,j}(T_{1}(k,p),T_{2}(n,j)),%
\end{array}
\label{Li} \\
\text{for }i\text{, }p<i\leq n,\text{and all }j\in\{p,\ldots,i-1\}\}  \notag
\end{gather}
Note that $R_{i}$ and $L_{i}$ are the tasks with the largest index in $%
\{p,\ldots,i-1\}$ that satisfies the inequalities in (\ref{Ri}) and (\ref{Li}%
) respectively. It is clear that $p\leq R_{i}<i,$ $p\leq L_{i}<i.$

A special case of (\ref{Ri}) and (\ref{Li}) arises when $p$ is the first
task of a BP $\{k,\ldots ,n\}$, i.e., $p=k.$ Then, according to the
definitions above, we obtain the following inequalities, which will be used
in our later results: 
\begin{gather}
\sigma _{k,R_{i}}(a_{k},d_{R_{i}})\leq \sigma _{k,m}(a_{k},d_{m}),
\label{(Ri)_example} \\
\text{for }i\text{, }k<i\leq n,\text{ and all }m\in \{k,\ldots ,i-1\}  \notag
\end{gather}%
\begin{gather}
\sigma _{k,L_{i}}(a_{k},T_{2}(n,L_{i}))\geq \sigma _{k,m}(a_{k},T_{2}(n,m)),
\label{(Li)_example} \\
\text{for }i\text{, }k<i\leq n,\text{and all }m\in \{k,\ldots ,i-1\}.  \notag
\end{gather}

After introducing the above definitions and notations, we are now ready to
introduce three important lemmas, which will be used to prove our main
theorem.

\begin{lem}
\label{No_left_critical_tasks}Let tasks $\{k,\ldots ,n\}$ form a BP on the
optimal sample path of \textbf{P1} and task $r>k$ be the first
right-critical task in $\{k,\ldots ,n\}$. If $\sigma _{k,r}(a_{k},d_{r})\geq
\sigma _{k,L_{r}}(a_{k},a_{L_{r}+1}),$ then there is no left-critical task
in $\{k,\ldots ,r-1\}.$
\end{lem}

{\color{black} }

\begin{lem}
\label{No_right_critical_tasks}Let tasks $\{k,\ldots ,n\}$ form a BP on the
optimal sample path of \textbf{P1}. Consider task $R_{i},$ for $i$, $k<i\leq
n.$ If $\sigma _{k,j}(a_{k},d_{j})\geq \sigma _{k,L_{j}}(a_{k},a_{L_{j}+1})$
and $\sigma _{k,j}(a_{k},a_{j+1})\leq \sigma _{k,R_{j}}(a_{k},d_{R_{j}}),$
for all $j$, $k<j<i$, then there is no right-critical task before task $%
R_{i}.$
\end{lem}

\begin{lem}
\label{Right_critical}Let tasks $\{k,\ldots ,n\}$ form a BP on the optimal
sample path of \textbf{P1}. If $\sigma _{k,i}(a_{k},a_{i+1})>\sigma
_{k,R_{i}}(a_{k},d_{R_{i}}),$ $\sigma _{k,j}(a_{k},d_{j})\geq \sigma
_{k,L_{j}}(a_{k},$ $a_{L_{j}+1})$ and $\sigma _{k,j}(a_{k},a_{j+1})\leq
\sigma _{k,R_{j}}(a_{k},d_{R_{j}}),$ for $i$, $k<i\leq n$, and for all $j$, $%
k<j<i$, then $R_{i}$ is right-critical.
\end{lem}

Before we introduce the main theorem, we would like to first summarize the
above three lemmas.

Lemma \ref{No_left_critical_tasks} provides the conditions under which there
are no left-critical tasks before the first right-critical task $r$ in a BP.

Lemma \ref{No_right_critical_tasks} provides the conditions under which
there are no right-critical tasks before a given task $R_{i}$ in a BP.

Lemma \ref{Right_critical} provides the conditions under which task $R_{i}$
in a BP is right-critical.

With the help of the above auxiliary results, we are able to establish the
following theorem, which can identify the first critical task in a BP on the
optimal sample path of \textbf{P1}:

\begin{thm}
\label{BP_First_critical}Let tasks $\{k,\ldots,n\}$ form a BP on the optimal
sample path of \textbf{P1}.

\emph{(i)} If%
\begin{equation}
\sigma_{k,j}(a_{k},d_{j})\geq\sigma_{k,L_{j}}(a_{k},a_{L_{j}+1}),
\label{Thm_1}
\end{equation}%
\begin{equation}
\sigma_{k,j}(a_{k},a_{j+1})\leq\sigma_{k,R_{j}}(a_{k},d_{R_{j}}),\text{ and }
\label{Thm_2}
\end{equation}%
\begin{equation}
\sigma_{k,i}(a_{k},a_{i+1})>\sigma_{k,R_{i}}(a_{k},d_{R_{i}}),  \label{Thm_3}
\end{equation}
for $i,$ $k<i\leq n$, and all $j$, $k<j<i$, then $R_{i}$ is the first
critical task in $\{k,\ldots,n\}$, and it is right-critical.

\emph{(ii)} If 
\begin{gather*}
\sigma_{k,j}(a_{k},d_{j})\geq\sigma_{k,L_{j}}(a_{k},a_{L_{j}+1}), \\
\sigma_{k,j}(a_{k},a_{j+1})\leq\sigma_{k,R_{j}}(a_{k},d_{R_{j}}),\text{ and}
\\
\sigma_{k,i}(a_{k},d_{i})<\sigma_{k,L_{i}}(a_{k},a_{L_{i}+1}),
\end{gather*}
for $i$, $k<i\leq n$, and all $j$, $k<j<i$, then $L_{i}$ is the first
critical task in $\{k,\ldots,n\}$, and it is left-critical.
\end{thm}

Let us look at the first part of Theorem \ref{BP_First_critical} again. The
first right-critical task of a BP on the optimal sample path of \textbf{P1}
can be correctly identified if we can find $i$ and $R_{i}$ which satisfy (%
\ref{Thm_1})-(\ref{Thm_3}). In essence, (\ref{Thm_1}) guarantees that there
is no left-critical task before $R_{i}$, (\ref{Thm_2}) guarantees that there
is no right-critical task before $R_{i},$ and (\ref{Thm_3}) guarantees that $%
R_{i}$ is a right-critical task. A similar argument applies to the second
part of the theorem. Therefore, the conditions in Theorem \ref%
{BP_First_critical} are not only sufficient but also necessary for
identifying the first critical task.

After obtaining the first critical task, either left-critical or
right-critical, the rest of the BP, can be considered as a new BP. Invoking
Lemma \ref{Critical_Condition}, the new BP starts at either the first
critical task's deadline (if it is right-critical) or the arrival time of
the next task after the first critical task (if it is left-critical).
Applying Theorem \ref{BP_First_critical} on the next BP, we are able to
identify its first critical task, which is the second critical task of the
original BP. Iteratively applying Theorem \ref{BP_First_critical} helps us
find \emph{all} critical tasks on the original optimal sample path. This
leads directly to an efficient algorithm which can identify all critical
tasks in BP $\{k,\ldots ,n\}$ on the optimal sample path of \textbf{P1}.
Meanwhile, as we have illustrated in Fig. \ref{Fig_Optimal_Structure}, after
identifying all critical tasks in BP $\{k,\ldots ,n\}$ on the optimal sample
path of \textbf{P1}, we can find all segments in $\{k,\ldots ,n\}$ with the
same energy function derivatives. Solving a $NE$ problem for each segment
and combining the solutions gives us the optimal solution to $Q(k,n)$.

The \emph{Generalized Critical Task Decomposition Algorithm} (GCTDA) which
identifies critical tasks and solves $Q(k,n)$ is as follows: \ 

\begin{tabular}{l}
\textbf{step\ 1 }$p=k;$%
\end{tabular}

\begin{tabular}{l}
\textbf{step 2} $i=p+1$, Solve $NE(p,p;T_{1}(k,p),T_{2}(n,p))$%
\end{tabular}

\begin{tabular}{l}
and $NE(p,p;T_{1}(k,p),d_{p});$%
\end{tabular}

\begin{tabular}{l}
\ \ Identify the first critical task in $(p,n)$%
\end{tabular}

\begin{tabular}{l}
\ \ while $(i\leq n)$%
\end{tabular}

\begin{tabular}{l}
$\ \{$Solve $NE(p,i;T_{1}(k,p),T_{2}(n,i))$%
\end{tabular}

\begin{tabular}{l}
\ \ and $NE(p,i;T_{1}(k,p),d_{i});$%
\end{tabular}

\begin{tabular}{l}
\ \ Compute $R_{i};$%
\end{tabular}

\begin{tabular}{l}
\ \ if $(\sigma_{p,i}(T_{1}(k,p),T_{2}(n,i))>%
\sigma_{p,R_{i}}(T_{1}(k,p),d_{R_{i}})$%
\end{tabular}

\begin{tabular}{l}
\ \ \ \{$R_{i}$ is the first right-critical task in $(p,n)$;%
\end{tabular}

\begin{tabular}{l}
$\ \ \ \ \tau_{j}^{\ast}=\tau_{j}(T_{1}(k,p),d_{R_{i}}),$ $j=p,\ldots,R_{i};$%
\end{tabular}

\begin{tabular}{l}
$\ \ \ \ x_{R_{i}^{\ast}}=d_{R_{i}};$%
\end{tabular}

\begin{tabular}{l}
$\ \ \ \ a_{j}=d_{R_{i}},$ for all $j$, s.t., $j>R_{i},a_{j}<d_{R_{i}};$%
\end{tabular}

\begin{tabular}{l}
$\ \ \ \ p=R_{i}+1;$ go to \textbf{step 2};\}%
\end{tabular}

\begin{tabular}{l}
\ \ Compute $L_{i};$%
\end{tabular}

\begin{tabular}{l}
\ \ if $(\sigma_{p,i}(T_{1}(k,p),d_{i})<%
\sigma_{p,L_{i}}(T_{1}(k,p),a_{L_{i}+1})$%
\end{tabular}

\begin{tabular}{l}
\ \ \ \ \{$L_{i}$ is the first left-critical task in $(p,n)$;%
\end{tabular}

\begin{tabular}{l}
$\ \ \ \ \tau_{_{j}}^{\ast}=\tau_{j}(T_{1}(k,p),a_{L_{i}+1}),$ $j=p,\ldots
,L_{i};$%
\end{tabular}

\begin{tabular}{l}
$\ \ \ \ x_{L_{i}}^{\ast}=a_{L_{i}+1};$%
\end{tabular}

\begin{tabular}{l}
$\ \ \ \ p=L_{i}+1;$ go to \textbf{step 2};\}%
\end{tabular}

\begin{tabular}{l}
$\ \ i=i+1;$%
\end{tabular}

\begin{tabular}{l}
\ \}%
\end{tabular}

\begin{tabular}{l}
$\ \ \tau_{j}^{\ast}=\tau_{j}(T_{1}(k,p),d_{n}),$ $j=p,\ldots,n;$%
\end{tabular}

\begin{tabular}{l}
\textbf{END}%
\end{tabular}

Note that GCTDA\ finds the critical tasks in a BP on the optimal sample path
of \textbf{P1} iteratively. The optimal departure times of these critical
tasks can be easily obtained using the results in Lemma \ref%
{Critical_Condition}. Finally, the optimal solution to the off-line problem: 
$\tau _{j}^{\ast },$ $j=1,\ldots N,$ is also calculated in GCTDA.

Regarding the complexity of our algorithm, the most time consuming part is
solving $NE(i,j;t_{1},t_{2})$. In the worst case, the optimal sample path is
a single BP containing $N-1$ critical tasks and the GCTDA algorithm may need
to solve $NE(i,j;t_{1},t_{2})$ $2N_{r}\ $times to identify each critical
task, where $N_{r}$ is the number of tasks remaining. Therefore, the worst
case complexity of the GCTDA algorithm is $O(N^{2})$.

\subsection{%
	Maximum Power Constraint%
}

In \textbf{P1}, we omitted the constraint: $\tau _{i}\geq \tau _{i\_\min },$
which is essentially the maximum transmission rate or transmission power
constraint for task $i$. This constraint is very important in real-world
scenarios because a transmitter simply cannot transmit above the maximum
transmission rate/power. We now formulate \textbf{P2}:
\begin{equation*}
\begin{tabular}{ll}
\textbf{P2:} & $\underset{\tau _{1}^{\prime },\ldots ,\tau _{N}^{\prime }}{%
	\min }\sum_{i=1}^{N}v_{i}\omega _{i}(\tau _{i}^{\prime })$ \\ 
$s.t.$ & $x_{i}^{\prime }=\max (x_{i-1}^{\prime },a_{i})+v_{i}\tau
_{i}^{\prime }\leq d_{i},$ $i=1,\ldots N$ \\ 
& $\tau _{i}^{\prime }>\tau _{i\_\min },\text{ }x_{0}=0.$%
\end{tabular}%
\end{equation*}%
Notice that the only difference between \textbf{P1} and \textbf{P2} is the
constraint on the control$.$ We use $\tau _{i}^{\prime \ast }$ and $%
x_{i}^{\prime \ast }$ to denote the optimal control and optimal departure
time of task $i$ in \textbf{P2}$,$ respectively. \\
It is easy to show that Lemmas \ref{DepartureIsDeadline} and \ref%
{DepartureLargerThanNextArrival} also apply to \textbf{P2}. Similar to how
we handled \textbf{P1}, we only need to consider a single BP $\{k,\ldots
,n\} $ in the optimal sample path of \textbf{P2}. We formulate the following
problem for BP $\{k,\ldots ,n\}$:%
\begin{equation*}
\begin{tabular}{ll}
$Q^{\prime }(k,n):$ & $\underset{\tau _{k}^{\prime },\ldots ,\tau
	_{n}^{\prime }}{\min }\sum_{i=k}^{n}v_{i}\omega _{i}(\tau _{i}^{\prime })$
\\ 
$\ \ \ s.t.$ & $x_{i}^{\prime }=a_{k}+\sum_{j=k}^{i}v_{i}\tau _{i}^{\prime
}\leq d_{i},$ $i=k,\ldots ,n$ \\ 
& $\tau _{i}^{\prime }\geq \tau _{i\_\min },$ $i=k,\ldots ,n$ \\ 
& $x_{i}^{\prime }\geq a_{i+1},\text{ }i=k,\ldots ,n-1.$%
\end{tabular}%
\end{equation*}%

In order to establish the connection between Problems \textbf{P1} and \textbf{P2}, we now introduce the following assumption, which will be used to derive the results in this subsection. 
Justifications for this assumption in transmission scheduling are provided in the appendix. 
\begin{assumption}
\label{A_tao_derivative}\emph{a)} If $\tau _{i\_\min }<\tau _{j\_\min },$
then $\dot{\omega}_{i}(\tau _{i\_\min })<\dot{\omega}_{j}(\tau _{j\_\min })$
and $\dot{\omega}_{i}(\tau )>\dot{\omega}_{j}(\tau );$ \emph{b)} If $\tau
_{i\_\min }\geq \tau _{j\_\min },$ then $\dot{\omega}_{i}(\tau _{i\_\min
})\geq \dot{\omega}_{j}(\tau _{j\_\min })$ and $\dot{\omega}_{i}(\tau )\leq 
\dot{\omega}_{j}(\tau ).$
\end{assumption}

Let 
\begin{equation*}
\tau _{\min }=\underset{i=k,\ldots ,n}{\inf }\tau _{i\_\min }.
\end{equation*}%

Under Assumption \ref{A_tao_derivative}, we introduce the following auxiliary lemma:

\begin{lem}
\label{connectionLemma}If $\exists \tau _{i}^{\ast }<\tau _{\min }$, then $%
Q^{\prime }(k,n)$ is infeasible.
\end{lem}

Lemma \ref{connectionLemma} establishes certain connections between the
power unconstrained problem $Q(k,n)$ and the power constrained problem $%
Q^{\prime }(k,n)$. When the problem is homogeneous, i.e., the cost functions
are identical among the tasks, we can easily derive that if $Q^{\prime
}(k,n) $ is feasible, then the optimal solution to $Q(k,n)$ must also yield
the maximum power constraint. In the inhomogeneous case, however, it is
possible that $Q(k,n)$ may return an optimal solution above the maximum
power constraint while $Q^{\prime }(k,n)$ is indeed feasible. When this
occurs, the optimal solution to $Q^{\prime }(k,n)$ would be close to $\tau
_{i\_\min },$ and the controller could simply apply $\tau _{i\_\min }$ as
the control since there is not much benefit to do optimization in this case.%

\subsection{Off-line Performance Comparisons}

Next, we will test the off-line performance of the GCTDA algorithm. In this
case, all task information including arrival times, deadlines, and number of
bits is known. For comparison purposes, we obtain numerical results for the
following algorithms:

\textbf{GCTDA}: Off-line algorithm knowing all task information exactly and
having full computational capability to solve $NE(i,j;t_{1},t_{2})$.

\textbf{GCTDA\_TL}: Off-line algorithm knowing all exact task information
and using pre-established tables to find an approximate solution to the
nonlinear algebraic system $NE(i,j;t_{1},t_{2}).$ The purpose of this
algorithm is to reduce the computational overhead associated with solving $%
NE(i,j;t_{1},t_{2}),$ at the cost of more energy consumption. Specifically,
we pre-calculate the derivatives of 1000 $\tau $ values for each energy
function $\omega _{i}(\tau )$ and save these data into tables. Using these
tables and binary search, we find approximate solutions to $%
NE(i,j;t_{1},t_{2})$ in GCTDA.

\textbf{MoveRight}: The algorithm proposed in \cite{GamNaPraUyZaInf02}. It is an iterative algorithm that converges to the optimal solution. We choose it for performance comparison purposes because to the best of our knowledge, it is the only other algorithm available for solving problems with task-dependent cost functions. 

In each experiment, in order to make the comparison fair, we use the same
setting (i.e., same arrival times, deadlines, task sizes, and energy
functions) for each algorithm. Note that what the ``best'' function solves
in the MoveRight algorithm is actually a nonlinear system $%
NE(i,i+1;t_{1},t_{2}) $. All experiments are done using a $1.8$GHz Athlon XP
processor.

The setting of the first experiment in Table I is as follows: $500$ tasks of
Poisson arrivals with mean inter-arrival time $5$s, each task has its own
deadlines (uniformly distributed between $[a_{i}+5,a_{i}+20]$ for task $i$),
task sizes are different, and the energy functions are the same. The GCTDA\
algorithm outperforms the MoveRight algorithm in terms of CPU time by two
orders of magnitude. Because the optimal sample path is likely to contain
multiple BPs and the energy functions are identical, GCTDA is very fast. We
terminated the MoveRight algorithm after $10000$ passes. It can be seen that
the MoveRight algorithm did not converge at this point yet (the cost is
still higher than the optimal cost returned by GCTDA.) Another observation
is that the solution of GCTDA\_TL is a good approximation to the one of
GCTDA. This makes GCTDA\_TL a good candidate for on-line control. However,
it 
\color{black}%
can be seen that GCTDA\_TL takes longer than GCTDA when the energy functions
are identical. The reason is that in this case, the nonlinear system $%
NE(i,j;t_{1},t_{2})$ becomes a linear system, which can be easily solved. So
there is no benefit in using the table lookup approximation approach.
However, when the energy functions are different, as we will see later, the
approximation method does help.

\begin{table}
\small
\centering
\resizebox{0.3\textwidth}{!}{
	\scalebox{0.7}{
	\tabcolsep=0.11cm
	\begin{tabular}{|c|c|c|}
	\hline
	& CPU time (sec) & Cost \\ 
	\hline 
	GCTDA & $0.031$ & $3.41919$ \\ 
	GCTDA\_TL & $1.579$ & $3.56494$ \\
	MoveRight & $54.704$ & $3.43264$ \\
	\hline
	\end{tabular}%
	}
}
\caption{\small Different task deadlines and identical energy
	functions}
\end{table}

\begin{table}
\small
\centering
\resizebox{.3\textwidth}{!}{
	\begin{tabular}{|c|c|c|}
	\hline
	& CPU time (sec) & Cost \\ \hline
	GCTDA & $12.516$ & $8.87826$ \\ 
	GCTDA\_TL & $2.469$ & $8.98774$ \\ 
	MoveRight & $200.703$ & $9.08324$ \\ 	\hline
	\end{tabular}%
}
\caption{\small Different deadlines and different energy functions}
\end{table}

\begin{table}
\centering
\small
\resizebox{.3\textwidth}{!}{
	\begin{tabular}{|c|c|c|}
	\hline
	& CPU time (sec) & Cost \\ \hline
	GCTDA & $0.593$ & $0.0687325$ \\ 
	GCTDA\_TL & $98.469$ & $0.0688239$ \\
	MoveRight & $61.969$ & $0.0837155$ \\ 	\hline
	\end{tabular}%
}
\caption{\small Identical deadlines and identical energy
	functions}
\end{table}

\begin{table}
\centering
\small
\resizebox{.3\textwidth}{!}{
	\begin{tabular}{|c|c|c|}
	\hline
	& CPU time (sec) & Cost \\ \hline
	GCTDA & $48$ & $0.1997$ \\ 
	GCTDA\_TL & $11.594$ & $0.200008$ \\
	MoveRight & $434.687$ & $0.72739$ \\	\hline
	\end{tabular}%
}
\caption{\small Identical deadlines and different energy
	functions}
\end{table}

In the next experiment for $500$ tasks in Table II, we keep the same setting
as above, except that we make the energy functions different for each task.
We terminate the MoveRight algorithm after $100$ passes. It can be seen that
in this experiment, GCTDA\_TL takes much less CPU time than GCTDA. Both of
them are much faster (by an order of magnitude) and MoveRight has not yet
converged.

In Table III, we make all $500$ tasks have the same deadline and the same
energy function. In this case, the optimal sample path contains a single BP.
We terminate the MoveRight algorithm after $10000$ passes. It can be seen
that at the time of termination, it was still far from converging to the
optimal solution. Again, the CPU time of GCTDA\_TL is higher than GCTDA,
since the energy functions are identical.

In Table IV, the setting is the same as above, except that we now consider $%
100$ tasks with different energy functions. We terminate the MoveRight
algorithm after $1000$ passes.

\section{On-line Controller Design}

We proved that our off-line algorithm GCTDA can return each critical task on
the optimal sample path correctly. Therefore, using it, we can get the
off-line optimal solution. We are also interested in designing good \emph{%
	on-line} controllers, in which case there are two difficulties: 1) lack of
future task 
\color{black}%
information; 2) high computational complexity in solving the nonlinear
equations.

To overcome the first difficulty, we design a Receding Horizon (RH)
controller assuming that at each decision point, the controller always has
some task information within a given RH window, and nothing beyond this
window. The size of the RH window $H$ can be measured either by time units
or the number of tasks. In this paper, we use the latter to measure the RH
window $H$. This RH\ window, together with the task information within it,
is often referred to as the \emph{planning horizon}. In contrast to the
planning horizon, the RH controller will apply controls over an \emph{action
	horizon, }which contains a subset of tasks over the planning horizon. Such
controllers have been proposed and analyzed in \cite{cgc-RH-stochastic} and 
\cite{MiaoCasRHTAC} for the homogeneous case that the cost functions are identical. In this
paper, we consider the RH control for the inhomogeneous case that the cost
functions are \emph{task-dependent}. 

As we will see later, the off-line results we obtained in previous sections provide insight to RH on-line controller design and performance evaluation. We now introduce some notations similar to the ones in \cite{MiaoCasRHTAC}.
Let $\tilde{x}_{t}$ be the departure time of task $t$ on the RH state
trajectory, which is also a decision point when the RH controller is invoked
with lookahead window $H$. Let $\tilde{\tau}_{t}$ be the control associated
with task $t$ as determined by the RH controller. When task $t+1$ starts a
new BP (i.e., $a_{t+1}>\tilde{x}_{t}$), then the RH controller does not need
to act until $a_{t+1}$ rather than $\tilde{x}_{t}$; for notational
simplicity, we will still use $\tilde{x}_{t}$ to represent the decision
point for task $t+1$ (i.e., the time when the control $\tilde{\tau}_{t+1}$
is determined). Let $h$ denote the last task included in the window that
starts at the current decision point $\tilde{x}_{t}$, i.e., 
\begin{equation*}
h=\arg \max {}_{r\geq t}\{a_{r}:a_{r}\leq \tilde{x}_{t}+H\}.
\end{equation*}%
Note that although the value of $h$ depends on $t,$ for notational
simplicity, we will omit this dependence and only write $h_{t}$ when it is
necessary to indicate dependence on $t$. When the RH controller is invoked
at $\tilde{x}_{t}$, it is called upon to determine $\tilde{\tau}_{i}$, the
control associated with task $i$ for all $i=t+1,\ldots ,h$, and let $\tilde{x%
}_{i}$ denote the corresponding departure time of task $i$\emph{\ }which is
given by $\tilde{x}_{i}=\max (\tilde{x}_{i-1},a_{i})+\tilde{\tau}_{i}v_{i}$.
The values of $\tilde{x}_{i}$ and $\tilde{\tau}_{i}$ are initially
undefined, and are updated at each decision point $\tilde{x}_{t}$ for all $%
i=t+1,\ldots ,h$. Control is applied to task $t+1$ only. That control and
the corresponding departure time are the ones showing in the final RH sample
path. In other words, for any given task $i$, $\tilde{x}_{i}$ and $\tilde{%
	\tau}_{i}$ may vary over different planning horizons, since optimization is
performed based on different available information. It is only when task $i$
is the next one at some decision point that its control and departure time
become final.

Given these definitions, we are now ready to discuss the worst case
estimation process to be used. If $h=N$, then the optimization process is
finalized, so we will only consider the more interesting case when $h<N$.
Then, our worst case estimation pertains to the characteristics of task $h+1$%
, the first one beyond the current planning horizon determined by $h$, i.e.,
its arrival time, deadline, and number of bits which are unknown. We define
task arrival times and task deadlines for $i=t+1,\ldots ,h+1$ as follows:%
\begin{eqnarray}
\tilde{a}_{i} &=&\left\{ 
\begin{array}{cl}
a_{i}, & \text{if }t+1\leq i\leq h \\ 
\tilde{x}_{t}+H, & \text{if }i=h+1%
\end{array}%
\right. \text{ }  \label{new_a} \\
\tilde{d}_{i} &=&\left\{ 
\begin{array}{cl}
d_{i}, & \text{if }t+1\leq i\leq h \\ 
\tilde{a}_{h+1}+\tau _{i\_\min }v_{h+1}, & \text{if }i=h+1%
\end{array}%
\right.  \label{new_d}
\end{eqnarray}%
In (\ref{new_a}), the arrival times of tasks $i=t+1,\ldots ,h$ are known and
we introduce a \textquotedblleft worst case\textquotedblright\ estimate for
the first unknown task beyond $\tilde{x}_{t}+H$, i.e., we set it to be the
earliest it could possibly occur. In (\ref{new_d}), the deadlines of tasks $%
i=t+1,\ldots ,h$ are known and we introduce a \textquotedblleft worst
case\textquotedblright\ estimate for the first unknown task's deadline to be
the tightest possible, since $\tau _{i\_\min }$ is the minimum feasible time
per bit. Note that $v_{h+1}$ is in fact unknown at time $\tilde{x}_{t}$, but
we will see that this does not affect our optimization process as the value
of $\tilde{d}_{h+1}$ is not actually required for analysis purposes. We
point it out that we do not have to worry about estimates for the unknown
tasks beyond $h+1$ (this is because of the FCFS nature of our system).

Therefore, the optimization problem the RH controller faces at time $\tilde{x%
}_{t}$ is over tasks $t+1,\ldots ,h$\emph{\ }with the added constraint that
they must all be completed by time $\tilde{a}_{h+1}=\tilde{x}_{t}+H$. This
is equivalent to redefining $\tilde{d}_{i}$ as%
\begin{equation}
\tilde{d}_{i}=\left\{ 
\begin{array}{cl}
d_{i}, & \text{if }t+1\leq i\leq h \\ 
\min (d_{h},\tilde{a}_{h+1}), & \text{if }i=h%
\end{array}%
\right.  \label{first_d_tilde}
\end{equation}%
Our on-line RH control problem at decision point $\tilde{x}_{t}$ will be
denoted by $\tilde{Q}(t+1,h)$ and is formulated as follows:%
\begin{equation*}
\begin{tabular}{ll}
$\tilde{Q}(t+1,h):$ & $\underset{\tilde{\tau}_{t+1},\ldots ,\tilde{\tau}_{h}}%
{\min }\text{ }\sum_{i=t+1}^{h}v_{i}\omega _{i}(\tilde{\tau}_{i})$ \\ 
$\ \ \ \ \ \ \ \ $s.t. & $\tilde{\tau}_{i}\geq 0,\text{ }i=t+1,\ldots ,h.$
\\ 
& $\tilde{x}_{i}=\max (\tilde{x}_{i-1},a_{i})+\tilde{\tau}_{i}v_{i}\leq 
\tilde{d}_{i},$ $\tilde{x}_{t}$ known$.\text{ }$%
\end{tabular}%
\end{equation*}%
where $\tilde{d}_{i}$ is defined in (\ref{first_d_tilde}). We also formulate
the on-line RH control problem with the maximum power constraint:%
\begin{equation*}
\begin{tabular}{ll}
$\tilde{Q}^{\prime }(t+1,h):$ & $\underset{\tilde{\tau}_{t+1},\ldots ,\tilde{%
		\tau}_{h}}{\min }\text{ }\sum_{i=t+1}^{h}v_{i}\omega _{i}(\tilde{\tau}_{i})$
\\ 
$\ \ \ \ \ \ \ \ $s.t. & $\tilde{\tau}_{i}\geq \tau _{i\_\min },\text{ }%
i=t+1,\ldots ,h.$ \\ 
& $\tilde{x}_{i}=\max (\tilde{x}_{i-1},a_{i})+\tilde{\tau}_{i}v_{i}\leq 
\tilde{d}_{i},\text{ }\tilde{x}_{t}\text{ }$known$.$%
\end{tabular}%
\end{equation*}

Similar to the RH problem in \cite{MiaoCasRHTAC}, $\tilde{Q}(t+1,h)$ may not
be feasible even if the off-line problem is feasible. This is due to the
worst-case estimation. One way of relaxing the worst-case estimation is to
use $\hat{h}$ (defined below)$,$ instead of $h$ in $\tilde{Q}$ above. Let%
\begin{eqnarray*}
\hat{x}_{j} &=&\max (\hat{x}_{j-1},a_{j})+\tau _{j\_\min }v_{j},\ \  \\
\ \hat{x}_{t} &=&\tilde{x}_{t},\text{ }j=t+1,\ldots ,h \\
S &=&\{j:t+1\leq j<h,\text{ } \\
\hat{x}_{i} &\leq &\min (d_{i},a_{j+1})\text{ for all }i,\text{ }t+1\leq
i\leq j\} \\
\hat{h} &=&\left\{ 
\begin{array}{cc}
\sup \text{ }S, & \text{if }S\neq \varnothing , \\ 
\infty , & \text{otherwise}%
\end{array}%
\right.
\end{eqnarray*}%
We then define $\hat{d}_{i}:$%
\begin{equation}
\hat{d}_{j}=\left\{ 
\begin{array}{cc}
d_{j}, & j=t+1,\ldots ,\hat{h}-1, \\ 
\min (d_{j},\tilde{a}_{j+1}), & j=\hat{h}.%
\end{array}%
\right. ,  \label{second_d_tilde}
\end{equation}%
and formulate problem $\hat{Q}(t+1,\hat{h}):$%
\begin{equation*}
\begin{tabular}{ll}
$\hat{Q}(t+1,\hat{h}):$ & $\underset{\tilde{\tau}_{t+1},\ldots ,%
	\tilde{\tau}_{\hat{h}}}{\min }\text{ }\sum_{i=t+1}^{\hat{h}}v_{i}\omega
_{i}(\tilde{\tau}_{i})$ \\ 
$\ \ \ \ \ \ \ \ $s.t. & $\tilde{\tau}_{i}\geq 0,\text{ }i=t+1,\ldots ,\hat{h%
}.$ \\ 
& $\tilde{x}_{i}=\max (\tilde{x}_{i-1},a_{i})+\tilde{\tau}_{i}v_{i}\leq 
\hat{d}_{i},$ $\tilde{x}_{t}$ known$.\text{ }$%
\end{tabular}%
\end{equation*}%
The RH control algorithm at each decision point $\tilde{x}_{t}$ is shown
in Table V. Note that $\tilde{Q}(t+1,h)$ and $\hat{Q}(t+1,\hat{h})$ essentially are smaller scale off-line optimization problems. This implies that at each on-line decision point, we shall use an off-line control algorithm, i.e., GCTDA or GCTDA\_TL, to solve $\tilde{Q}(t+1,h)$ and $\hat{Q}(t+1,\hat{h})$. 

In the next result, we discuss the feasibility of the proposed on-line RH control mechanism.  

\begin{thm}
\label{online_feasibility}If the off-line problem \textbf{P2} is feasible,
then the RH\ control in Table V is also feasible.
\end{thm}

Theorem \ref{online_feasibility} reveals that the RH control in Table V
guarantees feasibility when the off-line problem \textbf{P2} is feasible.
Next, we will analyze the performance of the proposed RH controller using
simulation. To overcome the high computational complexity, we use the
GCTDA\_TL algorithm, rather than GCTDA algorithm for on-line RH\ control. As
we have mentioned previously, the GCTDA\_TL algorithm uses some piecewise
constant functions to approximate the derivatives of the energy functions at
different $\tau $. Optimization can be approximated by searching efficiently
in a pre-established table containing these functions.

\begin{center}
\begin{tabular}{||l||l||}
\hline\hline
Step 1: & Solve $\tilde{Q}(t+1,h)$ and get $\widetilde{\tau }_{i}^{\ast },%
\tilde{x}_{i}^{\ast },$ \\ \hline\hline
& $i=t+1,\ldots ,h$ \\ \hline\hline
& If $\widetilde{\tau }_{i}^{\ast }\geq \tau _{i\_\min }$ for all $i$, \\ 
\hline\hline
& apply $\widetilde{\tau }_{t+1}^{\ast }$ to task $t+1$ and go to END. \\ 
\hline\hline
Step 2: & If $\hat{h}$ exists, solve $\hat{Q}(t+1,\hat{h})$ and
get $\widetilde{\tau }_{i}^{\ast },\tilde{x}_{i}^{\ast },$ \\ \hline\hline
& $i=t+1,\ldots ,\hat{h}.$ \\ \hline\hline
& If $\widetilde{\tau }_{i}^{\ast }\geq \tau _{i\_\min }$ for all $i$, \\ 
\hline\hline
& apply $\widetilde{\tau }_{t+1}^{\ast }$ to task $t+1$ and go to END. \\ 
\hline\hline
Step 3: & Apply $\tau _{t+1\_\min }$ to task $t+1$ \\ \hline\hline
END &  \\ \hline\hline
\end{tabular}

Table V:\ RH Control
\end{center}

In our simulation, all tasks have $512$ bytes and a fixed deadline, i.e., $%
d_{i}=a_{i}+d.$ The value of $d$ is set to $10s.$ In each figure below, we
run the experiment 1000 times, and 500 tasks are executed in each run. The
simulation is performed on a PC with a third generation Intel Core i5-3570K
Ivy Bridge 3.4GHz Quad-Core Desktop Processor. To quantify the deviation of
the RH cost from the optimal off-line cost, we define the cost difference
as: (RH cost - optimal off-line cost) / optimal off-line cost. The figures
below plot the average cost difference, worst-case cost difference,
best-case cost difference, and average calculation time versus the RH window
size $H$ in seconds.

\begin{figure}[tbph]
\begin{center}  
\includegraphics[width=8.5cm]{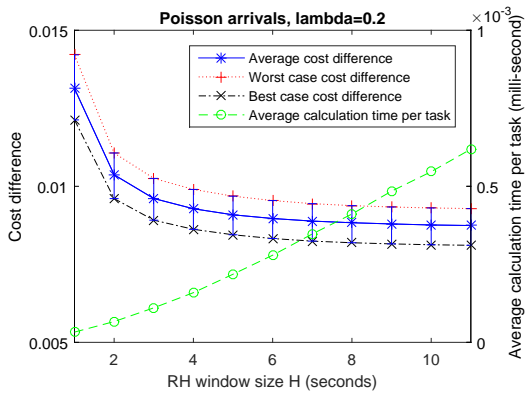}  
\vspace*{-10mm}
\caption{Simulation Results of Poisson Arrivals}.
\vspace*{-5mm} 
\label{online_simulation_Poisson}   
\end{center}  
\end{figure}

\begin{figure}[tbph]
\begin{center}  
\includegraphics[width=8.5cm]{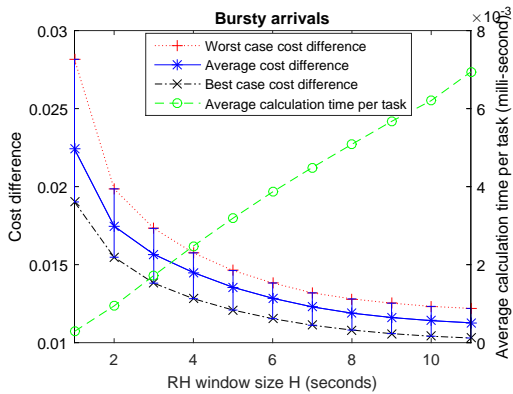}  
\vspace*{-10mm}
\caption{Simulation Results of Bursty Arrivals} \label{online_simulation_Bursty} 
\end{center}  
\end{figure}  

In Fig. \ref{online_simulation_Poisson}, we consider Poisson arrivals with $%
\lambda =0.2$. With RH window size $H$ varying from $1s$ to $11s$, the cost
differences are well below $2\%$. When $H$\ becomes larger, the cost
differences are reduced; the calculation time per task increases since at
each decision point, the optimization problem involves more tasks.

In the next experiment shown in Fig. \ref{online_simulation_Bursty}, we consider
bursty arrivals with the burst interval uniformly distributed over $[8,12]s$%
, the number of tasks in each burst chosen from $\{10,\ldots ,20\}$ with
equal probability, and the task intervals within the same burst uniformly
distributed within $[0,1]$. Although the cost differences are just slightly
higher than those in the Poisson case, the average calculation time per task
is now much larger. This is because in the bursty arrival case, a large
number of backlogged tasks are involved at each decision point.

Our simulation results show that the proposed RH control mechanism can not
only guarantee feasibility when the off-line problem is feasible, but also
achieve near optimal solutions.%

\section{Conclusions}

In this paper, we first study the \emph{Downlink Transmission Scheduling}
(DTS) problem. A simpler version of this problem has been studied in \cite%
{GamNaPraUyZaInf02} and \cite{TON-Energy-Efficient_Wireless}, where the
MoveRight algorithm is proposed. The MoveRight algorithm is an iterative
algorithm, and its rate of convergence is obtainable only when the cost
function is not \emph{task-dependent}. Compared with the work in \cite%
{TON-Energy-Efficient_Wireless} and \cite{GamNaPraUyZaInf02}, we deal with a
much harder problem: \emph{i)} our cost function is task-dependent and \emph{%
	ii)} each task has its own arrival time and deadline. This is essentially a
hard convex optimization problem with nondifferentiable constraints. By
analyzing the special structure of the optimal sample path, an efficient
algorithm, known as the Generalized Critical Task Decomposition\ Algorithm
(GCTDA), is proposed to solve the problem. Simulation results show that our
algorithm is more appropriate for real-time applications than the MoveRight
algorithm. Finally, we show that our results can be used in on-line control
to achieve near optimal solutions for Poisson and bursty arrivals.

\begin{center}
APPENDIX
\end{center}

\textbf{Proof of Lemma \ref{Convex_Transmission}:} Since $\omega _{i}(\tau )$
is strictly convex and differentiable,%
\begin{equation}
v_{1}\omega _{1}(\tau _{1})-v_{1}\omega _{1}(\tau _{1}^{^{\prime
	}})>v_{1}(\tau _{1}-\tau _{1}^{^{\prime }})\dot{\omega}_{1}(\tau
	_{1}^{^{\prime }})  \label{1}
	\end{equation}%
	\begin{equation}
	v_{2}\omega _{2}(\tau _{2})-v_{2}\omega _{2}(\tau _{2}^{^{\prime
		}})>v_{2}(\tau _{2}-\tau _{2}^{^{\prime }})\dot{\omega}_{2}(\tau
		_{2}^{^{\prime }})  \label{2}
		\end{equation}%
		Because $v_{1}\tau _{1}+v_{2}\tau _{2}=v_{1}\tau _{1}^{^{\prime }}+v_{2}\tau
		_{2}^{^{\prime }},$%
		\begin{equation}
		v_{1}(\tau _{1}-\tau _{1}^{^{\prime }})=-v_{2}(\tau _{2}-\tau _{2}^{^{\prime
			}})=C>0  \label{3}
			\end{equation}%
			Summing (\ref{1}) and (\ref{2}) above, and using (\ref{3}), we get:%
			\begin{gather*}
			v_{1}\omega _{1}(\tau _{1})+v_{2}\omega _{2}(\tau _{2})-v_{1}\omega
			_{1}(\tau _{1}^{^{\prime }})-v_{2}\omega _{2}(\tau _{2}^{^{\prime }})> \\
			C(\dot{\omega}_{1}(\tau _{1}^{^{\prime }})-\dot{\omega}_{2}(\tau
			_{2}^{^{\prime }}))
			\end{gather*}%
			Since $C>0$, and by assumption, $\dot{\omega}_{1}(\tau _{1}^{^{\prime }})-%
			\dot{\omega}_{2}(\tau _{2}^{^{\prime }})>0,$%
			\begin{equation*}
			v_{1}\omega _{1}(\tau _{1})+v_{2}\omega _{2}(\tau _{2})>v_{1}\omega
			_{1}(\tau _{1}^{^{\prime }})+v_{2}\omega _{2}(\tau _{2}^{^{\prime }}).\text{ 
			}{\small \blacksquare }
			\end{equation*}
			
			\textbf{Proof of Lemma \ref{Critical_Condition}: }We only prove part \emph{%
				(i)}. Part \emph{(ii)} can be proved similarly. Let $\tau _{k}^{^{\ast
				}},\ldots ,\tau _{n}^{\ast }$ be the optimal solution. By the definition of
				a left-critical task, we have $\dot{\omega}_{i}(\tau _{i}^{\ast })>\dot{%
					\omega}_{i+1}(\tau _{i+1}^{\ast }).$ Because tasks $i$ and $i+1$ are within
				a single BP, $x_{i}^{\ast }\geq a_{i+1}.$ Suppose $x_{i}^{\ast }>a_{i+1}.$
				Consider a feasible solution $\tau _{k}^{^{\prime }},\ldots ,\tau
				_{n}^{^{\prime }},$ $s.t.,$%
				\begin{align*}
				\tau _{j}^{^{\prime }}& =\tau _{j}^{^{\ast }},\text{ }j\neq i,\text{ }j\neq
				i+1, \\
				\tau _{i}^{^{\prime }}& <\tau _{i}^{\ast },\text{ }\tau _{i+1}^{^{\prime
					}}>\tau _{i+1}^{\ast },\text{ }\dot{\omega}_{i}(\tau _{i}^{^{\prime }})>\dot{%
					\omega}_{i+1}(\tau _{i+1}^{^{\prime }})
				\end{align*}%
				Note that such a feasible solution always exists as long as $\tau
				_{i}^{^{\prime }}$ and $\tau _{i+1}^{^{\prime }}$ are arbitrarily close to $%
				\tau _{i}^{\ast }$ and $\tau _{i+1}^{\ast }$ respectively. From Lemma \ref%
				{Convex_Transmission}, we get:%
				\begin{equation*}
				v_{i}\omega _{i}(\tau _{i}^{\ast })+v_{i+1}\omega _{i+1}(\tau _{i+1}^{\ast
				})>v_{i}\omega _{i}(\tau _{i}^{^{\prime }})+v_{i+1}\omega _{i}(\tau
				_{i+1}^{^{\prime }}).
				\end{equation*}%
				Since $\tau _{j}^{^{\prime }}=\tau _{j}^{^{\ast }},$ $j\neq i,$ $j\neq i+1,$
				using the above inequality, we get:%
				\begin{equation*}
				\sum_{i=k}^{n}v_{i}\omega _{i}(\tau _{i}^{\ast })>\sum_{i=k}^{n}v_{i}\omega
				_{i}(\tau _{i}^{^{\prime }}),
				\end{equation*}%
				which contradicts the assumption that $\tau _{k}^{^{\ast }},\ldots ,\tau
				_{n}^{\ast }$ is the optimal solution. Therefore, $x_{i}^{\ast }=a_{i+1}.$ $%
				{\small \blacksquare }$
				
				\textbf{Proof of Lemma \ref{NE_Properties}:} \emph{(i)} Since $\omega
				_{m}(\tau )$ is strictly convex, continuous and differentiable, $\dot{\omega}%
				_{m}(\tau )$ is also continuous. $NE(i,j;t_{1},t_{2})$ has a feasible
				solution, because $\dot{\omega}_{m}(\tau )$ is continuous and monotonically
				increasing, and $\tau $ can take any value in $(0,\infty ).$ Now, suppose
				there are two different solutions to $NE(i,j;t_{1},t_{2}):\tau
				_{i}(t_{1},t_{2}),\ldots ,\tau _{j}(t_{1},t_{2})$ and $\tau _{i}^{^{\prime
					}}(t_{1},t_{2}),\ldots ,\tau _{j}^{^{\prime }}(t_{1},t_{2}).$ Then, the
					common derivatives of these two solutions are different. Without loss of
					generality, we assume $\sigma _{i,j}(t_{1},t_{2})>\sigma _{i,j}^{^{\prime
						}}(t_{1},t_{2}),$ which means $\dot{\omega}_{m}(\tau _{m}(t_{1},t_{2}))>\dot{%
						\omega}_{m}(\tau _{m}^{^{\prime }}(t_{1},t_{2})),$ for any $m$, $i\leq m\leq
					j.$ By the convexity of $\omega _{m}(\tau )$, we obtain $\tau
					_{m}(t_{1},t_{2})>\tau _{m}^{^{\prime }}(t_{1},t_{2}),$ for any $m$, $i\leq
					m\leq j,$ and%
					\begin{align*}
					\sum_{m=i}^{j}\tau _{m}(t_{1},t_{2})v_{m}& =t_{2}-t_{1}> \\
					\sum_{m=i}^{j}\tau _{m}^{^{\prime }}(t_{1},t_{2})v_{m}& =t_{2}-t_{1},
					\end{align*}%
					which is a contradiction. Therefore, $NE(i,j;t_{1},t_{2})$ has a unique
					solution.
					
					\emph{(ii)} Let $\sigma_{i,j}(t_{1},t_{2})$, $\sigma_{i,j}(t_{3},t_{4})$ be
					the common derivative of $NE(i,j;$ $t_{1},t_{2})$ and $NE(i,j;t_{3},t_{4})$
					respectively, $0<\Delta=t_{2}-t_{1}<\Delta^{^{\prime}}=t_{4}-t_{3}.$ Let $%
					\tau_{i}(t_{1},t_{2}),\ldots,\tau_{j}(t_{1},t_{2})$, $\tau_{i}(t_{3},t_{4}),%
					\ldots,\tau_{j}(t_{3},t_{4})$ be the solution to $NE(i,j;t_{1},t_{2})$ and $%
					NE(i,j;t_{3},t_{4})$ respectively. We need to show $%
					\sigma_{i,j}(t_{1},t_{2})<$ $\sigma_{i,j}(t_{3},t_{4}).$ Suppose $%
					\sigma_{i,j}(t_{1},t_{2})\geq\sigma_{i,j}(t_{3},t_{4}).$ Then, by
					definition, we get $\dot{\omega}_{m}(\tau_{m}(t_{1},t_{2}))\geq\dot{\omega}%
					_{m}(\tau_{m}(t_{3},t_{4})),$ for any $m$, $i\leq m\leq j.$ By the convexity
					of $\omega _{m}(\tau)$, $\tau_{m}(t_{1},t_{2})\geq\tau_{m}(t_{3},t_{4}),$ $%
					i\leq m\leq j.$ Therefore, 
					\begin{align*}
					\sum_{m=i}^{j}\tau_{m}(t_{1},t_{2})v_{m} & =t_{2}-t_{1}=\Delta\geq \\
					\sum_{m=i}^{j}\tau_{m}(t_{3},t_{4})v_{m} & =t_{4}-t_{3}=\Delta^{^{\prime}}
					\end{align*}
					which contradicts our assumption that $\Delta<\Delta^{^{\prime}}.$
					Therefore, the common derivative of $NE(i,j;t_{1},t_{2})$ is a monotonically
					increasing function of $\Delta=t_{2}-t_{1}.$
					
					\emph{(iii)} It can be easily checked that $\tau_{i}(t_{1},t_{2}),\ldots,$ $%
					\tau_{p}(t_{1},t_{2})$ and $\tau_{p+1}(t_{1},$ $t_{2}),\ldots,%
					\tau_{j}(t_{1}, $ $t_{2})$ are the unique solutions to $%
					NE(i,p;t_{1},t_{1}+S_{ip})$ and $NE(p+1,j;t_{1}+S_{ip},t_{2})$ respectively$%
					. $ Therefore, by the definition of $NE(i,j;t_{1},t_{2})$ and $%
					\sigma_{i,j}(t_{1},t_{2}),$ we have $\sigma
					_{i,p}(t_{1},t_{1}+S_{ip})=\sigma_{i,j}(t_{1},t_{2})$ and $\sigma
					_{p+1,j}(t_{1}+S_{ip},t_{2})=\sigma_{i,j}(t_{1},t_{2})$, which implies $%
					\sigma_{i,p}(t_{1},t_{1}+S_{ip})=\sigma_{p+1,j}(t_{1}+S_{ip},t_{2})=%
					\sigma_{i,j}(t_{1},t_{2}).$
					
					\emph{(iv)} Let $S_{ip}$ be the same as in \emph{(iii)}. By assumption, $%
					c_{q}\neq c_{r}$ $\forall q,r\in \{1,2,3\},q\neq r.$ This implies that $%
					t_{3}\neq S_{ip},$ otherwise $c_{q}=c_{r}.$ From the monotonicity of $\sigma
					_{i,j}(t_{1},t_{2})$ shown in part \emph{(ii)}, we get $\min
					(c_{1},c_{2})<c_{3}<\max (c_{1},c_{2}).$\emph{\ }${\small \blacksquare }$
					
					\textbf{Proof of Lemma \textbf{\ref{No_left_critical_tasks}}:} Invoking
					Lemma \ref{Critical_Condition}, we get $x_{r}^{\ast }=d_{r}.$ Suppose there
					are left-critical tasks in $\{k,\ldots ,r-1\}$ and the closest left-critical
					task to $r$ is task $l$, $k\leq l<r.$ Invoking Lemma \ref{Critical_Condition}%
					, $x_{l}^{\ast }=a_{l+1}.$ By assumption, 
					\begin{equation}
					\sigma _{k,L_{r}}(a_{k},a_{L_{r}+1})\leq \sigma _{k,r}(a_{k},d_{r}).
					\label{19}
					\end{equation}%
					Because $l<r,$ from (\ref{(Li)_example})$,$ 
					\begin{equation}
					\sigma _{k,l}(a_{k},a_{l+1})\leq \sigma _{k,L_{r}}(a_{k},a_{L_{r}+1}).
					\label{20}
					\end{equation}%
					From (\ref{19}) and (\ref{20}), the following must be true:%
					\begin{equation}
					\sigma _{k,l}(a_{k},a_{l+1})\leq \sigma _{k,r}(a_{k},d_{r})  \label{24}
					\end{equation}%
					When the equality holds in (\ref{24}), from \emph{(iii)} of Lemma \ref%
					{NE_Properties}, we get 
					\begin{equation}
					\sigma _{l+1,r}(a_{l+1},d_{r})=\sigma _{k,l}(a_{k},a_{l+1}).  \label{24.5}
					\end{equation}%
					When the inequality holds in (\ref{24}), from \emph{(iv)} of Lemma \ref%
					{NE_Properties}, we obtain%
					\begin{equation*}
					\sigma _{l+1,r}(a_{l+1},d_{r})>\sigma _{k,r}(a_{k},d_{r}).
					\end{equation*}%
					From (\ref{24}) and the above inequality, we get 
					\begin{equation}
					\sigma _{l+1,r}(a_{l+1},d_{r})>\sigma _{k,l}(a_{k},a_{l+1}).  \label{24.8}
					\end{equation}%
					Combining (\ref{24.5}) and (\ref{24.8}), we get%
					\begin{equation}
					\sigma _{l+1,r}(a_{l+1},d_{r})\geq \sigma _{k,l}(a_{k},a_{l+1})  \label{10}
					\end{equation}%
					Since there is no right-critical or left-critical task in $\{l+1,\ldots
					,r-1\}$, invoking Lemma \ref{corollary_rate_equal}, we get%
					\begin{align*}
					\dot{\omega}_{s}(\tau _{s}^{\ast })& =\dot{\omega}_{s+1}(\tau _{s+1}^{\ast
					})=\sigma _{l+1,r}(a_{l+1},d_{r}),\text{ } \\
					\forall s& \in \{l+1,\ldots ,r-1\}.
					\end{align*}%
					From the definition of left-critical tasks, we get 
					\begin{equation}
					\dot{\omega}_{l}(\tau _{l}^{\ast })>\dot{\omega}_{l+1}(\tau _{l+1}^{\ast
					})=\sigma _{l+1,r}(a_{l+1},d_{r}).  \label{15}
					\end{equation}%
					We consider two cases:
					
					\emph{Case 1:} $k=l$. Then, $\sigma_{k,l}(a_{k},a_{l+1})=\dot{\omega}%
					_{l}(\tau_{l}^{\ast}).$ Inequalities (\ref{10}) and (\ref{15}) contradict
					each other.
					
					\emph{Case 2:} $k<l$. We will use a contradiction argument to show that
					there must exist a right-critical task $m$, $k\leq m<l$, i.e., 
					\begin{equation*}
					\dot{\omega}_{m}(\tau_{m}^{\ast})<\dot{\omega}_{m+1}(\tau_{m+1}^{\ast}).
					\end{equation*}
					Suppose such a task $m$ does not exist. i.e., 
					\begin{equation}
					\dot{\omega}_{m}(\tau_{m}^{\ast})\geq\dot{\omega}_{m+1}(\tau_{m+1}^{\ast }),%
					\text{ }k\leq m<l.  \label{21}
					\end{equation}
					Let 
					\begin{equation*}
					y_{m}=\left\{ 
					\begin{array}{ll}
					a_{k}, & m=k-1 \\ 
					x_{m}^{\ast} & m=k,\ldots,l.%
					\end{array}
					\right.
					\end{equation*}
					Inequality (\ref{21}) is equivalent to the following:%
					\begin{align}
					\sigma_{m,m}(y_{m-1},y_{m}) & \geq\sigma_{m+1,m+1}(y_{m},y_{m+1}),\text{ }
					\label{22} \\
					\text{for }m & =k,\ldots,l-1.  \notag
					\end{align}
					We will use a recursive proof next:
					
					\emph{Step 1}: Letting $m=k$ in (\ref{22}), we have 
					\begin{equation}
					\sigma_{k,k}(y_{k-1},y_{k})\geq\sigma_{k+1,k+1}(y_{k},y_{k+1}).  \label{22.5}
					\end{equation}
					When the equality holds in (\ref{22.5}), invoking part \emph{(iii)} of Lemma %
					\ref{NE_Properties}, we get%
					\begin{equation}
					\sigma_{k,k+1}(y_{k-1},y_{k+1})=\sigma_{k+1,k+1}(y_{k},y_{k+1}).
					\label{22.51}
					\end{equation}
					When the inequality holds in (\ref{22.5}), invoking part \emph{(iv)} of
					Lemma \ref{NE_Properties}, we get%
					\begin{equation}
					\sigma_{k,k+1}(y_{k-1},y_{k+1})>\sigma_{k+1,k+1}(y_{k},y_{k+1}).
					\label{22.52}
					\end{equation}
					Combining (\ref{22.51}) and (\ref{22.52}) above, we have 
					\begin{equation}
					\sigma_{k,k+1}(y_{k-1},y_{k+1})\geq\sigma_{k+1,k+1}(y_{k},y_{k+1}).
					\label{22.7}
					\end{equation}
					
					\emph{Step 2:} Letting $m=k+1$ in (\ref{22}), we have%
					\begin{equation*}
					\sigma _{k+1,k+1}(y_{k},y_{k+1})\geq \sigma _{k+2,k+2}(y_{k+1},y_{k+2}).
					\end{equation*}%
					Combining (\ref{22.7}) and the above inequality, we obtain%
					\begin{equation*}
					\sigma _{k,k+1}(y_{k-1},y_{k+1})\geq \sigma _{k+2,k+2}(y_{k+1},y_{k+2}).
					\end{equation*}%
					Similarly to the derivation of (\ref{22.51}), (\ref{22.52}), and (\ref{22.7}%
					), we can get%
					\begin{equation*}
					\sigma _{k,k+2}(y_{k-1},y_{k+2})\geq \sigma _{k+2,k+2}(y_{k+1},y_{k+2}).
					\end{equation*}%
					Repeating the process up to step $l-k,$ we obtain%
					\begin{equation*}
					\sigma _{k,l}(y_{k-1},y_{l})\geq \sigma _{l,l}(y_{l-1},y_{l}).
					\end{equation*}%
					Since task $l$ is left-critical, from the definition of $y_{m}$ and Lemma %
					\ref{Critical_Condition}, $y_{l}=x_{l}^{\ast }=a_{l+1},$ and the above
					inequality is equivalent to%
					\begin{equation}
					\sigma _{k,l}(a_{k},a_{l+1})\geq \dot{\omega}_{l}(\tau _{l}^{\ast }).
					\label{23}
					\end{equation}%
					From (\ref{23}) and (\ref{15}), we get%
					\begin{equation*}
					\sigma _{k,l}(a_{k},a_{l+1})>\sigma _{l+1,r}(a_{l+1},d_{r}).
					\end{equation*}%
					Invoking \emph{(iv)} of Lemma \ref{NE_Properties}, we obtain 
					\begin{equation*}
					\sigma _{k,l}(a_{k},a_{l+1})>\sigma _{k,r}(a_{k},d_{r})
					\end{equation*}%
					which contradicts (\ref{24}). Therefore, there must exist a task $m$, $k\leq
					m<l$, s.t., 
					\begin{equation*}
					\dot{\omega}_{m}(\tau _{m}^{\ast })<\dot{\omega}_{m+1}(\tau _{m+1}^{\ast })
					\end{equation*}%
					By Definition \ref{Definition1}, task $m$ is a right-critical task, which
					contradicts our assumption that task $r$ is the \emph{first} right-critical
					task in $\{k,\ldots ,n\}.$ Therefore, there is no left-critical task before
					task $r.$ ${\small \blacksquare }$
					
					\textbf{Proof of Lemma \textbf{\ref{No_right_critical_tasks}:}} We use a
					contradiction argument to prove the lemma. Suppose there are right-critical
					tasks before task $R_{i}$ and the one with smallest index is task $r$, $%
					k\leq r<R_{i}.$ By assumption, $\sigma _{k,j}(a_{k},d_{j})\geq \sigma
					_{k,L_{j}}(a_{k},a_{L_{j}+1}),$ for all $j$, $k<j<i.$ When $r>k$, letting $%
					j=r$, we have $\sigma _{k,r}(a_{k},d_{r})\geq \sigma
					_{k,L_{r}}(a_{k},a_{L_{r}+1}).$ Then we can invoke Lemma \ref%
					{No_left_critical_tasks} to establish that there is no left-critical task in 
					$\{k,\ldots ,r-1\}$. Since there is also no right-critical task in $%
					\{k,\ldots ,r-1\},$ from Lemma \ref{corollary_rate_equal}, 
					\begin{equation}
					\dot{\omega}_{s}(\tau _{s}^{\ast })=\sigma _{k,r}(a_{k},d_{r}),\text{ }%
					\forall s\in \{k,\ldots ,r\}.  \label{7.57}
					\end{equation}%
					Since $k\leq r<R_{i},$ from (\ref{(Ri)_example})$,$ we have%
					\begin{equation}
					\sigma _{k,R_{i}}(a_{k},d_{R_{i}})\leq \sigma _{k,r}(a_{k},d_{r}).
					\label{7.5}
					\end{equation}%
					Then, from \emph{(ii)} of Lemma \ref{NE_Properties}, 
					\begin{equation}
					\sigma _{k,R_{i}}(a_{k},x_{R_{i}}^{\ast })\leq \sigma
					_{k,R_{i}}(a_{k},d_{R_{i}}).  \label{7.51}
					\end{equation}%
					Combining (\ref{7.5}) and (\ref{7.51}), we get%
					\begin{equation}
					\sigma _{k,R_{i}}(a_{k},x_{R_{i}}^{\ast })\leq \sigma _{k,r}(a_{k},d_{r}).
					\label{7.56}
					\end{equation}%
					Since $r$ is right-critical, from (\ref{7.57}) and Definition \ref%
					{Definition1}, 
					\begin{equation}
					\sigma _{k,r}(a_{k},d_{r})=\dot{\omega}_{r}(\tau _{r}^{\ast })<\dot{\omega}%
					_{r+1}(\tau _{r+1}^{\ast }).  \label{7.59}
					\end{equation}%
					From (\ref{7.56}) and (\ref{7.59}), there must exist at least one
					left-critical task in $\{r+1,\ldots ,R_{i}-1\};$ otherwise, from the
					definition of a left-critical task in Definition \ref{Definition1} and a
					simple contradiction argument, we have $\dot{\omega}_{s}(\tau _{s}^{\ast
					})\leq \dot{\omega}_{s+1}(\tau _{s+1}^{\ast }),\forall s\in \{r+1,\ldots
					,R_{i}-1\}.$ Using this result, (\ref{7.59}) and a similar method as in
					obtaining (\ref{23}), we can get $\sigma _{k,r}(a_{k},d_{r})<\sigma
					_{k,R_{i}}(a_{k},x_{R_{i}}^{\ast }),$ which contradicts (\ref{7.56}).
					
					Let the left-critical task with smallest index be $l.$ From Lemma \ref%
					{Critical_Condition}, $x_{l}^{\ast }=a_{l+1}.$ Similar to obtaining $\sigma
					_{k,r}(a_{k},d_{r})<\sigma _{k,R_{i}}(a_{k},x_{R_{i}}^{\ast })$ above, we
					can get 
					\begin{equation}
					\sigma _{k,r}(a_{k},d_{r})<\sigma _{k,l}(a_{k},x_{l}^{\ast })=\sigma
					_{k,l}(a_{k},a_{l+1}).  \label{7.60}
					\end{equation}%
					By assumption, $\sigma _{k,j}(a_{k},a_{j+1})\leq \sigma
					_{k,R_{j}}(a_{k},d_{R_{j}}),$ and setting $j=l$, we get%
					\begin{equation*}
					\sigma _{k,l}(a_{k},a_{l+1})\leq \sigma _{k,R_{l}}(a_{k},d_{R_{l}}).
					\end{equation*}%
					Since, from (\ref{(Ri)_example}), we have $\sigma
					_{k,R_{l}}(a_{k},d_{R_{l}})\leq \sigma _{k,r}(a_{k},d_{r}),$ combining this
					and the above inequality, we obtain%
					\begin{equation*}
					\sigma _{k,l}(a_{k},a_{l+1})\leq \sigma _{k,r}(a_{k},d_{r})
					\end{equation*}%
					which contradicts (\ref{7.60}) and completes the proof. ${\small %
						\blacksquare }$
					
					\textbf{Proof of Lemma \textbf{\ref{Right_critical}:}} Suppose task $R_{i}$
					is not right-critical. By assumption, 
					\begin{equation*}
					\sigma _{k,R_{i}}(a_{k},d_{R_{i}})<\sigma _{k,i}(a_{k},a_{i+1}),
					\end{equation*}%
					Because $x_{R_{i}}^{\ast }\leq d_{R_{i}},$ from \emph{(ii)} of Lemma \ref%
					{NE_Properties}, 
					\begin{equation}
					\sigma _{k,R_{i}}(a_{k},x_{R_{i}}^{\ast })\leq \sigma
					_{k,R_{i}}(a_{k},d_{R_{i}}).  \label{7}
					\end{equation}%
					From the two inequalities above, we obtain:%
					\begin{equation}
					\sigma _{k,R_{i}}(a_{k},x_{R_{i}}^{\ast })<\sigma _{k,i}(a_{k},a_{i+1}).
					\label{9}
					\end{equation}%
					Invoking Lemma \ref{No_right_critical_tasks}, there is no right-critical
					task before task $R_{i}.$ Next, we use a contradiction argument to show that
					there must exist at least one right-critical task in $\{R_{i},\ldots ,i-1\}$%
					. Suppose there is no right-critical task in $\{R_{i},\ldots ,i-1\}$.
					Because there is no right-critical task in $\{k,\ldots ,i-1\}$, by
					Definition \ref{Definition1}, we have 
					\begin{equation*}
					\dot{\omega}_{m}(\tau _{m}^{\ast })\geq \dot{\omega}_{m+1}(\tau _{m+1}^{\ast
					}),\text{ }m=k,\ldots ,i-1.
					\end{equation*}%
					Let 
					\begin{equation*}
					y_{m}=\left\{ 
					\begin{array}{ll}
					a_{k}, & m=k-1 \\ 
					x_{m}^{\ast }, & m=k,\ldots ,i-1 \\ 
					a_{i+1}, & m=i%
					\end{array}%
					\right.
					\end{equation*}%
					The above inequality can be rewritten as%
					\begin{gather*}
					\sigma _{m,m}(y_{m-1},y_{m})\geq \sigma _{m+1,m+1}(y_{m},y_{m+1}), \\
					\text{ }m=k,\ldots ,i-1.
					\end{gather*}%
					Similar to the way of obtaining (\ref{23}) in proving Lemma \ref%
					{No_left_critical_tasks}, we obtain%
					\begin{equation*}
					\sigma _{k,R_{i}}(a_{k},x_{R_{i}}^{\ast })\geq \sigma _{k,i}(a_{k},a_{i+1}),
					\end{equation*}%
					which contradicts (\ref{9}).
					
					We showed there must exist at least one right-critical task in $%
					\{R_{i},\ldots,i-1\}.$ From the initial contradiction assumption, $R_{i}$ is
					not right-critical. When $i=R_{i}+1,$ the contradiction proof is completed.
					Next, we consider the case when $i>R_{i}+1.$
					
					Let $r$ be the closest right-critical task to $R_{i}$ in $\{R_{i},\ldots
					,i-1\}$. Since there is no right-critical task in $\{k,\ldots ,r-1\},$ by
					Definition \ref{Definition1}, 
					\begin{equation*}
					\dot{\omega}_{m}(\tau _{m}^{\ast })\geq \dot{\omega}_{m+1}(\tau _{m+1}^{\ast
					}),\text{ }m=k,\ldots ,r-1.
					\end{equation*}%
					Let 
					\begin{equation*}
					y_{m}=\left\{ 
					\begin{array}{ll}
					a_{k}, & m=k-1 \\ 
					x_{m}^{\ast }, & m=k,\ldots ,r%
					\end{array}%
					\right.
					\end{equation*}%
					The above inequality can be rewritten as%
					\begin{gather*}
					\sigma _{m,m}(y_{m-1},y_{m})\geq \sigma _{m+1,m+1}(y_{m},y_{m+1}), \\
					\text{ }m=k,\ldots ,r-1.
					\end{gather*}%
					Similar to the way of obtaining (\ref{23}) in proving Lemma \ref%
					{No_left_critical_tasks}, we obtain%
					\begin{equation*}
					\sigma _{k,R_{i}}(a_{k},x_{R_{i}}^{\ast })\geq \sigma
					_{k,r}(a_{k},x_{r}^{\ast })=\sigma _{k,r}(a_{k},d_{r}).
					\end{equation*}%
					From the above inequality and (\ref{7}), we obtain 
					\begin{equation*}
					\sigma _{k,r}(a_{k},d_{r})\leq \sigma _{k,R_{i}}(a_{k},d_{R_{i}})
					\end{equation*}%
					which contradicts the definition of $R_{i}$ in (\ref{Ri})$,$ since $%
					R_{i}<r<i.$ Therefore, task $R_{i}$ must be right-critical. ${\small %
						\blacksquare }$
					
					\textbf{Proof of Theorem \ref{BP_First_critical}:} We only prove part \emph{%
						(i)}. Part \emph{(ii)} can be proven similarly. The proof contains several
					steps: 1) using (\ref{Thm_1}), Lemma \ref{No_left_critical_tasks}, and
					setting $j=R_{i},$ we conclude that there is no left-critical task before $%
					R_{i};$ 2) using (\ref{Thm_1}), (\ref{Thm_2}), and Lemma \ref%
					{No_right_critical_tasks}, we establish that there is no right-critical task
					before $R_{i};$ 3) using (\ref{Thm_1}), (\ref{Thm_2}), (\ref{Thm_3}), and
					Lemma \ref{Right_critical}, it follows that $R_{i}$ is a right-critical
					task; and finally 4) combining the results established in the previous steps
					1)-3), we can obtain that $R_{i}$ is the first critical task in $\{k,\ldots
					,n\}$, and it is right-critical. ${\small \blacksquare }$
					
					\textbf{Justifications for Assumption} \ref{A_tao_derivative}: We only justify Part a).
					Part b) can be justified similarly. Using Shannon's theorem, $\tau _{i}$ can be represented by the following
					equation:%
					\begin{equation*}
					\tau _{i}=\frac{1}{B\log _{2}(1+\frac{s_{i}P}{N_{0}})}
					\end{equation*}%
					where $B$ is the bandwidth of the channel, $s_{i}$ is the \emph{%
						task-dependent} channel gain, $P$ is the transmission power, and $N_{0}$ is
					the power of the noise. Since $\frac{s_{i}P}{N_{0}}\gg 1$ in typical
					scenarios$,$ we can omit the $1$ above and represent $P$ in terms of $\tau
					_{i}:$%
					\begin{equation}
					P(\tau _{i})=\frac{N_{0}(2^{\frac{1}{B\tau _{i}}})}{s_{i}}  \label{p_of_tau2}
					\end{equation}%
					We assume that the maximum transmission power of each task is constant $P_{\max
					} $, and it determines $\tau _{i\_\min }$:%
					\begin{equation}
					P_{\max }=\frac{N_{0}(2^{\frac{1}{B\tau _{i\_\min }}})}{s_{i}}
					\label{P_max2}
					\end{equation}%
					Because 
					\begin{equation*}
					\omega _{i}(\tau _{i})=P(\tau _{i})\tau _{i},
					\end{equation*}%
					we use (\ref{p_of_tau2}) and (\ref{P_max2}) to get 
					\begin{equation}
					\dot{\omega}_{i}(\tau _{i})=P(\tau _{i})(1-\frac{1}{B\tau _{i}})=\frac{%
						N_{0}(2^{\frac{1}{B\tau _{i}}})}{s_{i}}(1-\frac{1}{B\tau _{i}})
					\label{w_derivative2}
					\end{equation}%
					\begin{equation}
					\dot{\omega}_{i}(\tau _{i})|_{\tau _{i}=\tau _{i\_\min }}=P_{\max }(1-\frac{1%
					}{B\tau _{i\_\min }})  \label{w_derivative_min2}
					\end{equation}
					
					Using (\ref{w_derivative_min2}) and $\tau _{i\_\min }<\tau _{j\_\min },$ we
					have 
					\begin{equation*}
					\dot{\omega}_{i}(\tau _{i\_\min })<\dot{\omega}_{j}(\tau _{j\_\min }).
					\end{equation*}
					
					Using (\ref{P_max2}) and $\tau _{i\_\min }<\tau _{j\_\min }$, we get%
					\begin{equation*}
					s_{i}>s_{j},
					\end{equation*}%
					i.e., the channel gain of task $i$ is greater than that of task $j$. Using (%
					\ref{w_derivative2}), we get%
					\begin{equation*}
					\dot{\omega}_{i}(\tau )>\dot{\omega}_{j}(\tau ).\text{ \ }{\small %
						\blacksquare }
					\end{equation*}
					
					\textbf{Proof of Lemma} \ref{connectionLemma}: Let us assume that there
					exists tasks $\{p,\ldots ,q\}$ $(k<p\leq q<n)$ in a BP $\{k,\ldots ,n\}$ of
					the optimal sample path of $Q(k,n)$ such that%
					\begin{eqnarray}
					\tau _{p-1}^{\ast } &\geq &\tau _{p-1\_\min },\text{ }\tau _{q+1}^{\ast
					}\geq \tau _{q+1\_\min },\text{ }  \label{tai_min_assumption} \\
					\tau _{i}^{\ast } &<&\tau _{\min },\text{ }i=p,\ldots ,q  \notag
					\end{eqnarray}%
					From Assumption \ref{A2}, 
					\begin{eqnarray}
					\dot{\omega}_{p-1}(\tau _{p-1}^{\ast }) &\geq &\dot{\omega}_{p-1}(\tau
					_{p-1\_\min }),\text{ }  \label{optimal_greater_than_min} \\
					\dot{\omega}_{q+1}(\tau _{q+1}^{\ast }) &\geq &\dot{\omega}_{q+1}(\tau
					_{q+1\_\min })  \notag
					\end{eqnarray}%
					Let 
					\begin{equation*}
					z=%
					\underset{i=k,\ldots ,n}{\arg \min }\text{ }\tau _{i\_\min }.
					\end{equation*}%
					
					Because $\tau _{i\_\min }\geq \tau _{z\_\min },$ we invoke Assumption \ref%
					{A_tao_derivative} and get%
					\begin{equation}
					\dot{\omega}_{i}(\tau _{\min })\leq \dot{\omega}_{z}(\tau _{\min }),\text{ }%
					i=p,\ldots ,q  \label{w_1}
					\end{equation}%
					From Assumption \ref{A2} and (\ref{tai_min_assumption}), we have%
					\begin{equation}
					\dot{\omega}_{i}(\tau _{i}^{\ast })<\dot{\omega}_{i}(\tau _{\min }),\text{ }%
					i=p,\ldots ,q  \label{w_2}
					\end{equation}%
					Combine (\ref{w_1}) and (\ref{w_2}) above, we have%
					\begin{equation}
					\dot{\omega}_{i}(\tau _{i}^{\ast })<\dot{\omega}_{z}(\tau _{\min }),\text{ }%
					i=p,\ldots ,q  \label{optimal_less_than_min}
					\end{equation}%
					Because $\tau _{p-1\_\min }\geq \tau _{\min }$, we invoke Assumption \ref%
					{A_tao_derivative} and get%
					\begin{equation}
					\dot{\omega}_{p-1}(\tau _{p-1\_\min })\geq \dot{\omega}_{z}(\tau _{\min })
					\label{p-1_larger}
					\end{equation}%
					Similarly, we use $\tau _{q+1\_\min }\geq \tau _{\min }$ and Assumption \ref%
					{A_tao_derivative} to get%
					\begin{equation}
					\dot{\omega}_{q+1}(\tau _{q+1\_\min })\geq \dot{\omega}_{z}(\tau _{\min })
					\label{q+1_larger}
					\end{equation}%
					Combining (\ref{optimal_greater_than_min}), (\ref{p-1_larger}), and (\ref%
					{q+1_larger}), we have%
					\begin{equation}
					\dot{\omega}_{p-1}(\tau _{p-1}^{\ast })\geq \dot{\omega}_{z}(\tau _{\min })%
					\text{ and }\dot{\omega}_{q+1}(\tau _{q+1}^{\ast })\geq \dot{\omega}%
					_{z}(\tau _{\min }).  \label{p_q_greater}
					\end{equation}%
					Then, we combine (\ref{optimal_less_than_min}) and (\ref{p_q_greater}) to get%
					\begin{eqnarray*}
						\dot{\omega}_{p-1}(\tau _{p-1}^{\ast }) &>&\dot{\omega}_{i}(\tau _{i}^{\ast
						})\text{ and }\dot{\omega}_{i}(\tau _{i}^{\ast })<\dot{\omega}_{q+1}(\tau
						_{q+1}^{\ast }),\text{ } \\
						i &=&p,\ldots ,q.
					\end{eqnarray*}%
					From the inequalities above, we have%
					\begin{equation*}
					\dot{\omega}_{p-1}(\tau _{p-1}^{\ast })>\dot{\omega}_{p}(\tau _{p}^{\ast })%
					\text{ and }\dot{\omega}_{q}(\tau _{q}^{\ast })<\dot{\omega}_{q+1}(\tau
					_{q+1}^{\ast }).
					\end{equation*}%
					Invoking Lemma \ref{Critical_Condition}, we have 
					\begin{equation}
					x_{p-1}^{\ast }=a_{p}\text{ and }x_{q}^{\ast }=d_{q}.  \label{start_and_end}
					\end{equation}%
					Because $x_{q}^{\ast }=x_{p-1}^{\ast }+\sum_{i=p}^{q}v_{i}\tau _{i}^{\ast },$
					we use (\ref{tai_min_assumption}) and (\ref{start_and_end}) to get:%
					\begin{eqnarray*}
						x_{q}^{\ast }-x_{p-1}^{\ast } &=&d_{q}-a_{p}=\sum_{i=p}^{q}v_{i}\tau
						_{i}^{\ast } \\
						&<&\sum_{i=p}^{q}v_{i}\tau _{\min }\leq \sum_{i=p}^{q}v_{i}\tau _{i\_\min }
					\end{eqnarray*}%
					For any feasible solution of $Q^{\prime }(k,n),$ we must have%
					\begin{equation*}
					\sum_{i=p}^{q}v_{i}\tau _{i}^{\prime }\leq d_{q}-a_{p}\leq
					\sum_{i=p}^{q}v_{i}\tau _{i\_\min },
					\end{equation*}%
					which implies that there exists at least one $\tau _{i}^{\prime },i=p,\ldots
					,q$, such that $\tau _{i}^{\prime }<\tau _{i\_\min }.$ This completes the
					proof that $Q^{\prime }(k,n)$ is infeasible. ${\small \blacksquare }$
					
					\textbf{Proof of Theorem} \ref{online_feasibility}: Because \textbf{P2} is
					feasible, we have 
					\begin{equation*}
					x_{i}^{^{\prime }\ast }\leq d_{i}.
					\end{equation*}%
					To prove the theorem, we only need to show that $\tilde{x}_{i}\leq
					x_{i}^{^{\prime }\ast }$ for $i=1,\ldots ,N.$ We use induction to prove it.
					
					1) When $t=0$, we are at the first decision point $\tilde{x}_{0}=a_{1}.$ We
					have two cases:
					
					Case 1.1:\ We apply $\tau _{1\_\min }$ to task 1. It is obvious that $\tilde{%
						x}_{1}\leq x_{i}^{^{\prime }\ast }.$
					
					Case 1.2:\ We apply control $\widetilde{\tau }_{1}^{\ast }$ obtained in
					either Step 1 or Step 2 of Table V to task 1. Without loss of generality,
					let us assume that $\widetilde{\tau }_{1}^{\ast }$ is obtained from Step 1
					of Table V. In this case, $\widetilde{\tau }_{i}^{\ast },i=1,\ldots ,h,$ are
					the solution to $\tilde{Q}(1,h)$ and $\widetilde{\tau }_{i}^{\ast }\geq \tau
					_{i\_\min }$ for all $i$. Therefore, $\widetilde{\tau }_{i}^{\ast }$ is also
					the solution to problem $\tilde{Q}^{\prime }(1,h).$ We have two subcases:
					
					Case 1.2.1: When the planning horizon contains the end of a busy period on
					the optimal sample path, it is trivial that $\widetilde{\tau }_{1}^{\ast
					}=\tau _{1}^{^{\prime }\ast }.$ Therefore, $\tilde{x}_{1}=x_{1}^{\prime \ast
				}\leq d_{1}.$
				
				Case 1.2.2: We now consider the more interesting case that $d_{i}\geq 
				\widetilde{a}_{i+1},i=1,\ldots ,h.$ To compare the RH problem and the
				off-line problem, we now add subscripts to indicate the starting and ending
				times of each problem$.$ In particular, we use $\tilde{Q}_{a_{1},a_{1}+H}^{%
					\prime }(1,h)$ to show that the starting transmission time of $\tilde{Q}%
				^{\prime }(1,h)$ is $a_{1}$ and the ending transmission time is $a_{1}+H.$
				Similarly, we use $Q_{a_{1},x_{h}^{\prime \ast }}^{\prime }(1,h)$ to show
				that the starting transmission time of $Q^{\prime }(1,h)$ is $a_{1}$ and the
				ending transmission time is $x_{h}^{^{\prime }\ast }.$ Because $%
				a_{h+1}>a_{1}+H$ and $d_{h}\geq a_{1}+H,$ we must have $x_{h}^{^{\prime
					}\ast }\geq a_{1}+H.$ Looking at $\tilde{Q}_{a_{1},a_{1}+H}^{\prime }(1,h)$
				and $Q_{a_{1},x_{h}^{\prime \ast }}^{\prime }(1,h),$ they are exactly the
				same, except that the ending transmission time of $Q_{a_{1},x_{h}^{\prime
						\ast }}^{\prime }(1,h)$ is potentially at a later time. Therefore, the
				optimal departure time of any task in $Q_{a_{1},x_{h}^{\prime \ast
					}}^{\prime }(1,h)$ must not be earlier than that in $\tilde{Q}%
					_{a_{1},a_{1}+H}^{\prime }(1,h),$ which means $\tilde{x}_{1}\leq
					x_{1}^{\prime \ast }.$
					
					2) Suppose that the RH\ controller is at decision point $\tilde{x}_{t}$, and 
					$\tilde{x}_{t}\leq x_{t}^{\prime \ast }.$ We also have two cases:
					
					Case 2.1:\ We apply $\tau _{t+1\min }$ to task $t+1$. It is obvious that $%
					\tilde{x}_{t+1}\leq x_{t+1}^{\prime \ast }.$
					
					Case 2.2: We apply control $\widetilde{\tau }_{t+1}^{\ast }$ obtained in
					either Step 1 or Step 2 of Table V to task 1. Without loss of generality,
					let us assume that $\widetilde{\tau }_{t+1}^{\ast }$ is obtained from Step 1
					of Table V. In this case, $\widetilde{\tau }_{i}^{\ast },i=t+1,\ldots ,h,$
					are the solution to $\tilde{Q}(t+1,h)$ and $\widetilde{\tau }_{i}^{\ast
					}\geq \tau _{i\_\min }$ for all $i$. Therefore, $\widetilde{\tau }_{i}^{\ast
				}$ is also the solution to problem $\tilde{Q}^{\prime }(t+1,h).$ We consider
				two subcases:
				
				Case 2.2.1: When the planning horizon contains the end of a busy period on
				the optimal sample path, i.e., there exists task $j\epsilon \{t+1,\ldots
				,h\} $ s.t. $d_{j}<\widetilde{a}_{j+1},$ we focus on tasks $\{t+1,\ldots
				,j\}.$ In this case, the controls of these tasks on the RH sample path are
				the solutions to problem $\tilde{Q}_{\tilde{x}_{t},d_{j}}^{\prime }(t+1,j)$,
				and the control of these tasks on the optimal sample path of \textbf{P2} are
				the solutions to problem $Q_{x_{t}^{\prime \ast },d_{j}}^{\prime }(t+1,j).$
				Looking at $\tilde{Q}_{\tilde{x}_{t},d_{j}}^{\prime }(t+1,j)$ and $%
				Q_{x_{t}^{\prime \ast },d_{j}}^{\prime }(t+1,j)$, they are identical, except
				that the starting transmission time of $\tilde{Q}_{\tilde{x}%
					_{t},d_{j}}^{\prime }(t+1,j)$ is potentially earlier than that of $%
				Q_{x_{t}^{\prime \ast },d_{j}}^{\prime }(t+1,j).$ Therefore, the optimal
				departure time of any task in $\tilde{Q}_{\tilde{x}_{t},d_{j}}^{\prime
				}(t+1,j)$ must be no later than that in $Q_{x_{t}^{\prime \ast
				},d_{j}}^{\prime }(t+1,j),$ which means $\tilde{x}_{t+1}\leq x_{t+1}^{\prime
				\ast }.$
			
			Case 2.2.2: $d_{i}\geq \widetilde{a}_{i+1},i=t+1,\ldots ,h.$ In this case,
			the controls of tasks $\{t+1,\ldots ,h\}$ on the RH sample path are the
			solutions to problem $\tilde{Q}_{\tilde{x}_{t},\tilde{x}_{t}+H}^{\prime
			}(t+1,h)$, and the control of these tasks on the optimal sample path of 
			\textbf{P2} are the solutions to problem $Q_{x_{t}^{\prime \ast
				},x_{h}^{\prime \ast }}^{\prime }(t+1,h).$ Because $a_{h+1}>\tilde{x}_{t}+H$
			and $d_{h}\geq \tilde{x}_{t}+H,$ we must have $x_{h}^{^{\prime }\ast }\geq 
			\tilde{x}_{t}+H.$ Looking at $\tilde{Q}_{\tilde{x}_{t},\tilde{x}%
				_{t}+H}^{\prime }(t+1,h)$ and $Q_{x_{t}^{\prime \ast },x_{h}^{\prime \ast
				}}^{\prime }(t+1,h),$ they are exactly the same, except that both the
				starting and ending transmission times of $\tilde{Q}_{\tilde{x}_{t},\tilde{x}%
					_{t}+H}^{\prime }(t+1,h)$ are potentially sooner than those of $%
				Q_{x_{t}^{\prime \ast },x_{h}^{\prime \ast }}^{\prime }(t+1,h)$. Therefore,
				the optimal departure time of any task in $\tilde{Q}_{\tilde{x}_{t},\tilde{x}%
					_{t}+H}^{\prime }(t+1,h)$ must be no later than that in $Q_{x_{t}^{\prime
						\ast },x_{h}^{\prime \ast }}^{\prime }(t+1,h),$ which means $\tilde{x}%
				_{t+1}\leq x_{t+1}^{\prime \ast }.$ ${\small \blacksquare }$


\end{document}